\documentclass[11pt]{article}

\usepackage{amsmath}
\usepackage{amssymb}
\usepackage{latexsym}
\usepackage{amsmath, amsfonts,amssymb, amsthm, euscript,makeidx,color,mathrsfs}
\usepackage{cite}

\newtheorem{assumption}{Assumption}
\def\qed{ \ \vrule width.2cm height.2cm depth0cm\smallskip}

\newcommand{\la}{\langle}
\newcommand{\ra}{\rangle}
\newcommand{\hP}{\hat\dbP}

\newcommand{\ba}{\begin{array}}
\newcommand{\ea}{\end{array}}
\newcommand{\be}{\begin{equation}}
\newcommand{\ee}{\end{equation}}
\newcommand{\bea}{\begin{eqnarray}}
\newcommand{\eea}{\end{eqnarray}}
\newcommand{\beaa}{\begin{eqnarray*}}
\newcommand{\eeaa}{\end{eqnarray*}}

\def\neg{\negthinspace}

\def\dbE{\mathbb{E}}
\def\dbF{\mathbb{F}}

\def\dbP{\mathbb{P}}
\def\dbR{\mathbb{R}}

\def\a{\alpha}

\def\g{\gamma}
\def\d{\delta}
\def\e{\varepsilon}

\def\m{\mu}
\def\n{\nu}
\def\si{\sigma}
\def\t{\tau}
\def\f{\varphi}
\def\th{\theta}
\def\o{\omega}

\def\G{\Gamma}
\def\D{\Delta}

\def\L{\Lambda}

\def\O{\Omega}
\def\cA{{\cal A}}

\def\cE{{\cal E}}
\def\cF{{\cal F}}
\def\cG{{\cal G}}

\def\hC{\mathbb{C}}
\def\hD{\mathbb{D}}
\def\hE{\mathbb{E}}
\def\hF{\mathbb{F}}
\def\hG{\mathbb{G}}
\def\hH{\mathbb{H}}

\def\hL{\mathbb{L}}

\def\hN{\mathbb{N}}

\def\hP{\mathbb{P}}
\def\hQ{\mathbb{Q}}
\def\hR{\mathbb{R}}
\def\hS{\mathbb{S}}

\def\sB{\mathscr{B}}
\def\sC{\mathscr{C}}

\def\sE{\mathscr{E}}

\def\sK{\mathscr{K}}

\def\sP{\mathscr{P}}

\def\sT{\mathscr{T}}

\def\sX{\mathscr{X}}

\def\no{\noindent}

\def\ss{\smallskip}
\def\ms{\medskip}
\def\bs{\bigskip}
\def\q{\quad}
\def\qq{\qquad}

\def\pa{\partial}
\def\cd{\cdot}
\def\cds{\cdots}
\def\lan{{\langle}}
\def\ran{{\rangle}}

\def\tr{\hbox{\rm tr}}

\def\qed{ \hfill \vrule width.25cm height.25cm depth0cm\smallskip}
\newcommand{\dfnn}{\stackrel{\triangle}{=}}
\newcommand{\basa}{\begin{assumption}}
\newcommand{\easa}{\end{assumption}}

\newcommand{\bas}{\begin{assum}}
\newcommand{\eas}{\end{assum}}

\def\liminf{\mathop{\underline{\rm lim}}}

\def\lan{\mathop{\langle}}
\def\ran{\mathop{\rangle}}

\def\essinf{\mathop{\rm essinf}}

\def\limP2{\,\mathop{\buildrel \Pi_2\over\longrightarrow\,}}

\def\pa{\partial}

\def\wt{\widetilde}
 \def\cd{\cdot}
\def\cds{\cdots}
\def\ae{\hbox{\rm -a.e.{ }}}
\def\as{\hbox{\rm -a.s.{ }}}

\def\tr{\hbox{\rm tr$\,$}}

\def\ind{{\perp\neg\neg\neg\perp}}

\def\dis{\displaystyle}
\def\wt{\widetilde}

\def\1{{\bf 1}}

\def\:{\!:\!}

at 9pt

\begin{document}

\newtheorem{thm}{Theorem}[section]
\newtheorem{lem}[thm]{Lemma}
\newtheorem{cor}[thm]{Corollary}
\newtheorem{prop}[thm]{Proposition}
\newtheorem{rem}[thm]{Remark}
\newtheorem{eg}[thm]{Example}
\newtheorem{defn}[thm]{Definition}
\newtheorem{assum}[thm]{Assumption}

\renewcommand {\theequation}{\arabic{section}.\arabic{equation}}
\def\thesection{\arabic{section}}

\title{\bf A General Conditional McKean-Vlasov Stochastic Differential Equation}
\author{Rainer Buckdahn$^{1,2}$,\,\, Juan Li$^{3}$,\,\, Jin Ma$^4$ \\
 {$^1$\small Laboratoire de Math\'{e}matiques de Bretagne Atlantique, Univ Brest,}\\
	{\small UMR CNRS 6205, 6 avenue Le Gorgeu, 29200 Brest, France.}\\
	{$^2$\small  School of Mathematics, Shandong University,}
	{\small Jinan 250100, P. R. China.}\\
	{$^3$\small  School of Mathematics and Statistics, Shandong University, Weihai,}
	{\small Weihai 264209, P. R. China.}\\
{$^4$\small  Department of
Mathematics, University of Southern California,}
	{\small Los Angeles, 90089, USA.}\\
 {\small{\it E-mails: rainer.buckdahn@univ-brest.fr,\,\ juanli@sdu.edu.cn,\,\ jinma@usc.edu.}}}
\renewcommand{\thefootnote}{\fnsymbol{footnote}}
\footnotetext[1]{Rainer Buckdahn is supported in part by the ``FMJH Program Gaspard Monge in optimization and operation research", and the ANR (Agence Nationale de la Recherche), France project ANR-16-CE40-0015-01. Juan Li is supported by the NSF of P.R. China (NOs. 12031009, 11871037), National Key R and D Program of China (NO. 2018YFA0703900). Jin Ma is supported in part by US NSF grants \#DMS-1908665.

\ \ $^{**}$Juan Li is the corresponding author.}
\date{August 07, 2021}
\maketitle

\begin{abstract}
In this paper we consider a class of {\it conditional McKean-Vlasov SDEs} (CMVSDE for short).
Such an SDE can be considered as an extended version of McKean-Vlasov SDEs with common noises, as well as the
general version of the so-called {\it conditional mean-field SDEs} (CMFSDE) studied previously by the authors \cite{BLM, MSZ},
but with some fundamental differences. In particular, due to the lack of compactness of the iterated conditional laws, the existing arguments of Schauder's fixed point theorem do not  seem to apply in this situation, and the heavy nonlinearity on the conditional laws caused by change of probability measure adds more technical subtleties. Under some structure assumptions on the coefficients
of the observation equation, we prove the well-posedness of solution in the weak sense along a more direct approach. Our result is  the first that deals with McKean-Vlasov type SDEs involving state-dependent conditional laws.
\end{abstract}



\section{Introduction}
\label{sect-Introduction}
\setcounter{equation}{0}

In this paper we are interested in the well-posedness of the following general form of {\it conditional McKean-Vlasov} stochastic
differential equations (SDEs), defined on a certain filtered probability space $(\O, \cF, \hP, \hF=\{\cF_t\})$:
\bea
\label{SDE1}
\left\{\ba{lll}
\dis dX_t= b(t, \cd, X_{\cd\wedge t}, Y_{\cd\wedge t}, \m^{X|Y}_{\cd\wedge t})dt+\sum_{i=1}^2\si_i(t, \cd, X_{\cd\wedge t}, Y_{\cd\wedge t}, \m^{X|Y}_{\cd\wedge t})dB^i_t, ~ &X_0=x;\ms\\
dY_t= h(t, \cd, X_{\cd\wedge t}, Y_{\cd\wedge t}, \m^{X|Y}_{\cd\wedge t})dt+\hat\si  dB^2_t, &Y_0=0,
\ea\right.
\eea
where $b, h, \si_1, \si_2$ are measurable functions defined on appropriate spaces, $\hat\si$ is a constant,  $(B^1,B^2)$ is an
$(\mathbb{F},\hP)$-Brownian motion,
and $\mu_t^{X|Y}(\cd):=\hP\{X_t\in \cdot~|\mathcal{F}_t^Y\}$ denotes the regular conditional distribution of $X_t$ given $\cF^Y_t=\si\{Y_s:
s\le t\}$.

Special forms of SDE (\ref{SDE1}) have appeared in many applications, especially those involving partial informations, and have been studied by the authors in different co-authorships in the past (see, for example, \cite{BLM}, \cite{MSZ}). In many of these applications
the conditional law appears in the form of conditional expectations $\hE[X_t|\cF^Y_t]$, in the spirit of the nonlinear filtering problems, and hence often refer to as {\it conditional mean-field} SDEs. Consequently, the coefficients of these SDEs depend either linearly
on  $\hE[X_t|\cF^Y_t]$ (see, e.g., \cite{MSZ}), or linearly on the law of  $\hE[X_t|\cF^Y_t]$  (see, e.g., \cite{BLM}). The SDEs with
a general coupling between the solutions and their conditional law in the coefficients such as (\ref{SDE1}) have not been completely explored yet in the literature. In fact, there seem to be some fundamental difficulties when the usual solution methods
are employed.

SDE (\ref{SDE1}) can also be viewed from another angle. Assuming for example $h\equiv 0$, then $\cF^Y\equiv \cF^{B}$,
and the SDE becomes the so-called McKean-Vlasov SDE with common noise.
We refer to \cite{BCEH, CDL, CG, HSS, LS} and the references cited therein for various recent studies for SDEs with similar natures and their applications. In this case, two facts are worth noting: 1) the probability measure determining the conditional law is fixed throughout; and 2) the conditioning filtration is given exogenously, and is independent of the state $X$.  The case when the coefficient $h\neq 0$, however, is quite different. Since the ``observation" process $Y$ depends on $X$, the conditioning filtration $\cF^Y$ becomes {\it state-dependent}, whence
endogenous. Among other complications caused by such a ``coupling" nature, one of the severe consequences is that the conditional laws $\{\m^{X|Y}\}$ are no longer compact, loosing an important technical basis of the well-posedness arguments for McKean-Vlasov SDEs with common noises (see, e.g., \cite{HSS}).

To illustrate this point, let us ask the following simple question one would encounter naturally in constructing any  iteration scheme
in seeking the solution for SDE (\ref{SDE1}):
Given a pair of random variables taking values in any metric space, does the strong convergence
$(X^n,Y^n)\to(X, Y)$ imply the convergence
$\m^{X^n|Y^n}\to \m^{X|Y}$, in the sense of probability distributions? The answer to this question unfortunately negative. For example, let
$X^n \equiv U$, where $U$ is a random variable such that Var$(U)>0$ (whence $\hP\{U\neq \hE^{\hP}[U]\}>0$), and $Y^n=\frac{1}{n}X_n=\frac{1}{n}U$. Then $X^n\to X=U$,  and $Y^n\to Y\equiv0$, as $n\to \infty$. 
Obviously, for suitable non-constant, bounded measurable function $f$, we have, for any $n$,
$\hE[f(X^n)|\mathcal{F}^{Y^n}]=f(U) \neq \hE[f(U)]= \hE[f(X)|\mathcal{F}^Y]$.
This shows, in particular, that $\m^{X^n|Y^n}$ does not converge to $\m^{X|Y}$. We note, however, that in the usual common noise case the conditioning $\si$-field
is fixed (i.e., $Y^n=Y=B^2$), so such a problem does not occur.

In light of the nonlinear filtering theory, a tempting remedy to ``fix" the conditioning fitration is to consider the so-called {\it reference
measure} $\hQ^0$, a {\it prior} probability measure that is equivalent to $\hP$ but under which $(B^1, Y)$ is a Brownian motion.
But doing so would lead to another dilemma: The conditional law in (\ref{SDE1}) is defined under the original probability $\hP$ (under
which $(B^1, B^2)$ is a Brownian motion), not the reference measure $\hQ^0$. The two conditional laws can be connected via the
Bayes rule (known as the {\it Kallianpur-Strieble} formula), but will inevitably cause
some serious technical issues,  especially when the conditional law $\m^{X|Y}$ (under  $\hP$) is now a part of the solution
 of the CMVSDE (\ref{SDE1}).

Our plan of attack is based on the following basic ideas. We shall design an iteration scheme which would include the conditional law $\m^{X|Y}$ as a component, considered as
{\it measure-valued} process defined on an appropriate space where
the weak convergence can be more conveniently analyzed.
More specifically, we shall argue that $\m^{X|Y}$ is a  measure-valued process that has continuous paths in the space of
probability laws under the Wasserstein metric, and that it can be identified as part of the fixed point, along with the processes
$(X, Y)$.
The main difficulty in implementing such an idea is that throughout the process we need to use the reference probability $\hQ^0$, via
the Kallianpur-Strieble formula. This leads to some new technicalities that are not commonly seen in the existing literature of nonlinear filtering or the McKean-Vlasov SDEs with common noises. In particular, it seems that a certain boundedness of the Girsanov kernel involved in connecting  reference measure and the original ones (in both directions) becomes inevitable, and it essentially amounts to asking for
a pathwise bound for the solution of a linear SDE (or a martingale), which is next to impossible. As a consequence, we shall
impose a structural assumption on the observation drift coefficient $h$, and we hope to be able to remove such restrictions in our future works.

This paper is organized as follows. In Section 2 we introduce the basic notations, definitions, and assumptions. In particular, we shall define the processes of the conditional laws, and establish some basic facts on its path regularities in terms of the Wasserstein metric.
In Section 3 we introduce our solution scheme and give some justifications of our main ideas. In Section 4 we establish our fundamental estimates.
 In the Sections 5 and 6 we prove the existence and uniqueness (in law) of the weak solution, respectively.

\section{Preliminaries}
\setcounter{equation}{0}

Throughout this paper we denote $\hC_T:= \hC([0,T], \hR)$, and   $\hP^0$ to be the Wiener measure on $\hC_T$.
We shall consider the following {\it canonical space} $(\O^0, \cF^0, \hQ^0)$:
\bea
\label{canonical}
\O^0:=\hC_T^2 := \hC([0,T];\hR^2),\q \cF^0:=\sB(\hC^2_T), \q \hQ^0:=\hP^0\otimes\hP^0.
\eea
In the above, $\sB(\hC^2_T)$ denotes the Borel $\si$-field on $\hC_T^2$. Furthermore, we denote $(B^1, Y)$ to be the canonical process, that is, $(B^1_t, Y_t)(\o)=(\o^1(t), \o^2(t))$, $t\in[0,T]$, where $\o=(\o^1, \o^2)\in \hC^2_T$. Then $(B^1, Y)$ is a 2-dimensional-Brownian motion under $\hQ^0$. Also, we define $\hF^0=\{\cF^0_t\}_{t\in[0,T]}:=\{\sB_t(\hC_T^2)\}_{t\in[0,T]}$, where $\sB_t(\hC^2_T):=
\si\{\o(\cd\wedge t): \o\in\hC^2_T\}$, to be the natural filtration generated by $(B^1, Y)$, and we denote
$\hF:=\overline{\hF^0}^{\hQ^0}$, the augmentation
of $\hF^0$ under $\hQ^0$, so that $\hF$ satisfies the {\it usual hypotheses}.

Now let $(\sX, d)$ be any metric space, and $\sB(\sX)$ the topological Borel $\si$-field on $\sX$.  For any sub-$\si$-field $\cG
\subseteq \cF^0$,  and $p\ge 1$, we denote  $\hL^p_{\cG}(\sX)$ to be the space of all random variables $\xi\mapsto \sX$,
such that $\xi$ is $\cG$-measurable and for any/some $e\in \sX$, $\hE^{\hQ^0}[d(e,\xi)^p]<\infty$. Similarly, for a sub-filtration $\hG\subseteq\hF$, and $p\ge 1$, we let
 $\hL^p_\hG([0,T];\sX)$ be the space of all $\sX$-valued, $\hL^p$-integrable, $\hG$-adapted processes
on $[0,T]$.  Furthermore, %
we denote
$\sC_T(\sX)$ to be all $\sX$-valued continuous functions defined on $[0,T]$, and
denote $\hL^0_\hG(\sC_T(\sX))$ to be the space of all $\sX$-valued, $\hG$-adapted
continuous processes. Finally, for any
$1\leq p<\infty$, we define
\bea
\label{Sp}
\left\{\ba{lll}
\hS_{\hG}^p(\sX):=\{z\in \hL^0_{\hG}(\sC_T(\sX)):
\hE^{\hQ^0}\Big[\sup\limits_{t\in[0,T]}d(z_t,e)^p\Big]<+\infty, ~\exists e\in \sX\}; \ms\\
\hS^{\infty-}_{\hG}(\sX):=\bigcap_{p\ge 1}\hS^p_{\hG}(\sX).
\ea\right.
\eea

Let us now denote $\sP(\sX)$ to be the space of all  probability measures on the metric space $\sX$, and $\sP_p(\sX)=\{\g\in\sP(\sX): \int_{\sX} d(z,e)^p \gamma(dz)<+\infty, \exists\ e\in \sX\}\subset \sP(\sX)$, $p\ge 1$. For $p\ge 1$, we endow $\sP_p(\sX)$ with the {\it p-Wasserstein metric}:
\bea
\label{Wass1}
W^p_p(\g_1,\g_2)&:=&\inf\{\int_{\sX^2}d(z_1,z_2)^p\rho(dz_1 dz_2) :\rho\in\sP(\sX^2),\rho(\cdot\times\sX)=\g_1,\rho(\sX\times\cdot)=\g_2\} \nonumber\ms\\
&=&\inf\{\hE^{\hQ^0}[d(\xi_1,\xi_2)^p]\ : \xi_1,\xi_2\in \hL^1_{\cF^0}(\sX)\}.
\eea

Recall that, if the metric space $(\sX,d)$ is complete, then also $(\sP_p(\sX), W_p)$ is complete, for all $p\geq 1$.

In what follows we shall focus on the case $p=1$. It is well-known that, for $\sX =\mathbb{R}$, $(\sP_1(\hR),W_1(\cdot,\cdot))$ is a complete and separable  metric space.
Furthermore, since $\O^0=\hC^2_T$ is Polish, we know that for each $t\in[0,T]$, the regular conditional probability $\hQ^{\o^2}_t(\,\cd\,):=\hQ^0(\,\cd\,|\cF^Y_t\}(\o^2)$ exists, that is, for any $A\in \cF^0$, $\o^2\mapsto \hQ^{\o^2}_t(A)$ is $\sB_t(\hC_T)/\sB(\hR)$ measurable and, for
any $\o^2\in\O^0$, $\hQ^{\o^2}_t(\cd)$ is a probability measure. Since $\sB_t(\hC_T)$ is generated by the paths $Y_{\cd\wedge t}(\o^2)=\o^2(\cd\wedge t)$, $\o^2\in\hC_T$, we will denote
$\hQ^{\o^2}_t=\hQ^{Y_{\cd\wedge t}}=\hQ^{\o^2_{\cd\wedge t}}$, when there is no confusion.

Now for any random variable $\xi$ defined on $(\O^0, \cF^0, \hQ^0)$, and $t\in[0,T]$, we consider the  regular conditional distribution:
\bea
\label{PY}
\hP^{Y_{\cd\wedge t}}_\xi(\, \cd \,)(\o^2)=\hQ^0[\xi\in \cd\, |\cF^Y_t](\o^2)=\hQ^{\o^2}_t\circ \xi^{-1}(\cd)=\hP^{\o^2_{\cd\wedge t}}_\xi(\cd)\in \sP_1 (=: \sP_1\left(\mathbb{R}\right)).
\eea
We would like to show that the mapping $(t, \o^2)\mapsto \hP_\xi^{\o^2_{\cd\wedge t}}
(\cd)\in \sP_1$ actually defines a measure-valued process as it should.
More precisely, we have the following result.
 \begin{lem}
\label{EYmeas}
Let $\xi$ be a random variable defined on $(\O^0, \cF^0, \hQ^0)$. Then for each $t\in[0,T]$, the mapping $\o^2\mapsto \hP^{Y_{\cd\wedge t}}_\xi(\cd)(\o^2)=\hP^{\o^2_{\cd\wedge t}}_\xi(\cd)$
is $\sB_t(\hC_T)/\sB(\sP_1)$-measurable.
\end{lem}

\noindent{\it Proof.} To begin with, note that $\sB(\sP_1)=\si\{B_r(\m): \m\in\sP_1, r>0\}$, where $B_r(\m)=\{\n\in\sP_1: W_1(\n,\m)\le r\}$. Thus,
$ [\hP^{\o^2_{\cd\wedge t}}_\xi]^{-1}(B_r(\m))=\{\o^2\in\hC_T: W_1(\hP^{\o^2_{\cd\wedge t}}_\xi, \m)\le r\}$.
Next, we recall the {\it Kantorovich-Rubinstein formula} (cf. \cite{KA} or \cite{KR}):
\bea
\label{W1}
W_1(\hP^{\o^2_{\cd\wedge t}}_\xi, \m)=\sup\Big\{\Big|\int_{\hR}\varphi d\hP^{\o^2_{\cd\wedge t}}_\xi-\int_{\hR}\varphi d{\mu}\Big|:~ \varphi\in \mbox{Lip}_1(\hR)\Big\},
\eea
where Lip$_1(\hR)$ is the space of all Lipschitz functions with Lipschitz constant 1.
We claim that there exists a countable subset $\L\subset$ Lip$_1(\hR)$ such that (\ref{W1}) can be replaced by
\bea
\label{W2}
W_1(\hP^{\o^2_{\cd\wedge t}}_\xi, \m)=
 \sup\Big\{\Big|\int_{\hR}\varphi d\hP^{\o^2_{\cd\wedge t}}_\xi-\int_{\hR}\varphi d{\mu}\Big|:~ \f\in \L\Big\},
\eea
and we can then conclude that
\beaa
\label{meas}
&&[\hP^{\o^2_{\cd\wedge t}}_\xi]^{-1}(B_r(\m))=\{\o^2\in\hC_T: W_1(\hP^{\o^2_{\cd\wedge t}}_\xi, \m)\le r\}\nonumber\\
&=& \{\o^2\in\hC_T: \sup_{\f\in \L}\Big|\int_{\hR}\varphi d\hP^{\o^2_{\cd\wedge t}}_\xi-\int_{\hR}\varphi d{\mu}\Big|\le r\}\nonumber\\
&=&\bigcap_{\f\in \L}\{\o^2\in\hC_T: \Big|\int_{\hR}\varphi d\hP^{\o^2_{\cd\wedge t}}_\xi-\int_{\hR}\varphi d{\mu}\Big|\le r\}\\
&= &\bigcap_{\f\in \L}\{\o^2\in\hC_T: \Big|\hE^{\hQ^0}[\f(\xi)|\cF^Y_{t}](\o^2)-\bar\f\Big|\le r\}= \bigcap_{\f\in \L} \hE^{\hQ^0}[\f(\xi)|\cF^Y_{t}]^{-1}(B_r(\bar\f)), \nonumber
\eeaa
where $\bar\f:=\int_{\hR}\varphi d{\mu}\in\hR$. Since
for each $\f\in \L$, the mapping $\o^2\mapsto \hE^{\hQ^0}[\f(\xi)|\cF^Y_{t}](\o^2)$ is $\sB_t(\hC_T)$-measurable, we conclude that $[\hP^{\o^2_{\cd\wedge t}}_\xi]^{-1}(B_r(\m)) \in \sB_t(\hC_T)$.

It remains to find the countable subset $\L\in \mbox{Lip}_1(\hR)$ so that (\ref{W2}) holds. To this end, we consider the following subset
of $\hC^1(\hR)$:
\bea
\label{H1}
\hH_0:=\{ f\in \hC^1(\hR): |f|_{\hH}^2 :=|f(0)|^2+\int_{\hR}|f'(y)|^2e^{-|y|}dy <\infty\},
\eea
and let $\hH=\overline{\hH}_0$, the closure of $\hH_0$ under the norm $|\cd|_{\hH}$. Then $(\hH, |\cd|_{\hH})$ is a separable Banach
space (in fact, a Hilbert space), and clearly Lip$_1(\hR)\subset \hH$.

Now, for each $f\in \hH$, define
$\f_f(x)=f(0)+\int_0^x [(f'(y)\wedge 1)\vee (-1)]dy$, $x\in\hR$, then $\f_f\in \mbox{Lip}_1(\hR)$. Furthermore, since $(\hH, |\cd|_{\hH})$ is separable, there exists a countable dense subset $\L_{\hH}\subset \hH$, and then it is not hard to check that
$\L:=\{\f_f:f\in \L_\hH\}$ is a countable dense subset of Lip$_1(\hR)$ under the norm $|\cd|_{\hH}$.  Consequently, for any $\m\in\sP_1(\hR)$ and $f\in \mbox{Lip}_1(\hR)$, we can find $\{f_n\}\subset \L\subset \mbox{Lip}_1(\hR)$ such that $|f-f_n|_{\hH}\to 0$, as
$n\to \infty$. Since $f$ and $f_n$'s are all absolutely continuous, we have
\bea
\label{appro}
\Big|\int_\hR fd\m-\int_\hR f_nd\m\Big|\le |f(0)-f_n(0)|+\Big|\int_\hR\Big|\int_0^x|f'(y)-f_n'(y)|dy\Big|\m(dx)\Big|.
\eea
Note that for each $x\in\hR$, we have $\Big[\int_0^x|f'(y)-f_n'(y)|dy\Big]^2\le |x|e^{|x|}|f-f_n|_{\hH}\to 0$, as $n\to \infty$, and
since $\{f, f_n, n\in\hN\}\subset \mbox{Lip}_1(\hR)$, we have $\big|\int_0^x|f'(y)-f_n'(y)|dy\big|\le 2|x|$. But $\m\in\sP_1(\hR)$ implies
that $\int_{\hR} |x|\m(dx)<\infty$. The Dominated Convergence Theorem and  (\ref{appro}) thus imply that $\int_{\hR}f_n d\m
\to \int_{\hR}fd\m$, as $n\to\infty$. This, together with (\ref{W1}), easily leads to (\ref{W2}). The proof is now complete.
 \qed

We remark that Lemma \ref{EYmeas} does not imply directly that the mapping $(t,\o^2)\mapsto \hP^{Y_{\cd\wedge t}}_\xi$ is
$\sB([0,T])\otimes \sB(\hC_T)/\sB(\sP_1)$  jointly measurable. But as we shall argue in the next section (see, also \cite{BLM}), that for fixed $\o^2\in
\hC_T$, the mapping
$t\mapsto \hP^{\o^2_{\cd\wedge t}}_\xi$ is a $\sP_1$-valued continuous function, which then renders the desired joint measurability.
Throughout our paper we shall focus on the case $\sX=\sC_T(\sP_1)$, the space of all $\sP_1$-valued continuous functions defined
on $[0,T]$.  The set $\hS^p_{\hF^Y}(\sC_T(\sP_1))$ defined by (\ref{Sp}) as well as the set
$\hS^{\infty-}_{\hF^Y}(\sC_T(\sP_1))$ will therefore be particularly useful in our discussion.

To conclude this section we introduce the following standard assumptions on the coefficients $b, \si_1, \si_2$, and $h$ of SDE (\ref{SDE1}). For convenience, in what follows we shall assume $\si_1=\si$, $\si_2=0$, and $\hat\si= 1$.
\begin{assum}
\label{assump1}
The function $\f=(b, \si, h):\left[0,T\right]\times\Omega\times\mathbb{C}_T^2\times\sC_T(\sP_1)\mapsto \hR^3$ is bounded, progressively measurable, and  for some constant $C>0$, it holds
that
\bea
\label{Lip}
&&|\f(t, x, y, \m)-\f(t, x', y, \m')|
\le C\Big[\sup_{s\in[0,t]}|x_s-x'_s| +\sup_{s\in[0,t]} W_1(\m_s, \m'_s)\Big], \\
&&\qq\qq\qq\qq\qq\qq\qq t\in[0,T], ~x, x', y\in \hC_T, ~\m,\m'\in\sC_T(\sP_1).\nonumber
\eea
\end{assum}

We also recall the notion of a {\it weak solution} to
SDE (\ref{SDE1}), which will be the main objective of this paper.
\begin{defn}
\label{weaksol}
A six-tuple  $(\O,\cF,\hF,\hP,(B^1,B^2),(X,Y))$ is called  a weak solution of (\ref{SDE1}), if:

{\rm (i)}  $(\Omega,\mathcal{F},\mathbb{F},P)$ is a filtered probability space satisfying the usual assumptions;

\ss
{\rm (ii)} $(B^1,B^2)$ is an $(\hF, \hP)$-Brownian motion;

\ss
{\rm (iii)} $(X,Y)\in \hL^2_{\mathbb{F}}([0,T];\hR^2)$ such that all terms in (\ref{SDE1}) are well-defined and (\ref{SDE1}) holds
for all  $t\in[0,T]$, $\hP$-a.s.
\end{defn}

\begin{rem}
\label{remark0}
{\rm It is worth noting that Definition \ref{weaksol} only defines processes $(X, Y)$, along with a probability set-up including the Brownian motion $(B^1, B^2)$. The conditional law $\m^{X|Y}$ then comes naturally as the function of $X$ and $Y$, under probability
$\hP$. But the example in the introduction shows that, unless some more structural information on the process $\m^{X|Y}$ is known,
the simple minded iteration scheme will likely fail. The main idea of our solution scheme is to add the conditional law $\m^{X|Y}$ into the iteration process itself to help the convergence analysis.
\qed}
\end{rem}

%

\section{The Solution Scheme}
\setcounter{equation}{0}

In this section we introduce the iteration scheme that will lead to the desired weak solution. A key
element in this scheme is the process of conditional laws, $\m^{X|Y}=\{\m^{X|Y}_t\}$  which, by Lemma \ref{EYmeas}, is a $\sP_1(\hR)$-valued measurable process, and will be used to ``decouple" the SDEs for $X$ and $Y$ in (\ref{SDE1}).
In light of the analysis in our previous work \cite{BLM}, we shall argue that it is actually
a $\sP_1(\hR)$-valued continuous process. That is, $\m^{X|Y}\in \sC_T(\sP_1)$. We therefore shall start our scheme by
considering $\m$ as a free variable  taking values in $\hL^0_{\hF^Y}(\sC_T(\sP_1))$, and then try to find the desired conditional law
 by a fixed-point argument. All our arguments are essentially independent of the drift coefficient $b$, under Assumption \ref{assump1}. Thus, for notational simplicity, in what follows we shall assume that $b\equiv 0$, as adding it back does not cause substantial difficulties.

To begin with,  for $\m\in \hL^0_{\hF^Y}(\sC_T(\sP_1))$, we consider the following simplified system of SDEs on $(\O^0, \cF^0, \hQ^0)$:
\bea
\label{SDE2}
\left\{\ba{lll}
dX_t= \si(t,  X_{\cd\wedge t}, Y_{\cd\wedge t}, \m_{\cd\wedge t})dB^1_t, \q  &X_0=x;\ms\\
dL_t= h(t, X_{\cd\wedge t}, Y_{\cd\wedge t}, \m_{\cd\wedge t})L_tdY_t, &L_0=1.
\ea\right.
\eea

Since $\m$ is $\hF^Y$-adapted, we can write $\m_{\cd\wedge t}=\Phi_t(Y_{\cd\wedge t})$, $t\in[0,T]$, $\hQ^0$-a.s., for some progressively measurable functional $\Phi:[0,T]\times \hC_T\mapsto \sC_T(\sP_1)$. But
$Y$ is part of the canonical process, SDE (\ref{SDE2}) can be thought of as one that has random and functional
type coefficients. Thus under Assumption \ref{assump1}, it has a unique strong solution on the probability space $(\O^0, \cF^0, \hQ^0)$, and we denote it by $(X^\m, L^\m)$. Since  $h$ is bounded, we see that the process $L^\m$ is an
$(\hF, \hQ^0)$-martingale, and can be written as the Dol\'eans-Dade stochastic exponential:
\bea
\label{Dade}
L^\m_t=\exp\Big\{\int_0^t  h(s, X_{\cd\wedge s}, Y_{\cd\wedge s}, \m_{\cd\wedge s})dY_s-\frac12\int_0^t |h(s, X_{\cd\wedge s}, Y_{\cd\wedge s}, \m_{\cd\wedge s})|^2ds\Big\},
\eea
$t\in [0,T]$. Moreover, since also $\si$ is bounded, it is not hard to show that $(X^\m, L^\m)\in
\hS^{\infty-}_{\hF}(\hR^2)$. Furthermore, as a strong solution, there exists a measurable non-anticipating functional
$\Psi: \hC_T^2\times \sC(\sP_1)\mapsto \hC_T^2$ such that $(X^\m, L^\m)=\Psi(B^1, Y, \m)$, $\hQ^0$-a.s.

%
Next, the $\hQ^0$-martingale $L^\m$
defines a new probability measure on $(\O^0, \cF)$: $\hP^\m(d\o):= L^\m_T \hQ^0(d\o)$. Then, under the new probability $\hP^\m$, the process $(B^1, B^2=Y-\int_0^\cd h(s, X^\m_{\cd\wedge t}, Y_{\cd\wedge t},$ $\m_{\cd\wedge t})dt)$ is a Brownian motion. Furthermore,  we
denote the  {\it regular conditional probability distribution} of the process $X^\m$, given $\cF^Y$, under the probability measure
$\hP^\m$ by $\tilde\m_t$, $t\in[0,T]$. Since $\tilde \m$ is obviously uniquely determined for each $\m\in  \hL^0_{\hF^Y}(\sC_T(\sP_1))$, we
can then define the so-called {\it solution mapping} by $\sT(\m):= \tilde \m$. That is, $\sT$ is a mapping from $\hL^0_{\hF^Y}(\sC_T(\sP_1))$ to $ \hL^0_{\hF^Y}([0,T]; \sP_1)$, and by the {\it Kallianpur-Strieble formula} we see that, for $A\in \sB(\O)$ and
$t\in[0,T]$,  one has
\bea
\label{muth}
[\sT(\m)]_t(A)=\tilde\m_t(A)\dfnn \hP^\m\{X^\m_t\in A|\cF^Y_t\}=\frac{\hE^{\hQ^0}\big[L^\m_t{\bf 1}_{\{X^\m_t\in A\}} |\cF^Y_t\big]}
{\hE^{\hQ^0}\big[L^\m_t |\cF^Y_t\big]}.
\eea

Let us now assume that the mapping $\sT$ has  a fixed point. That is, there exists
$\hat \m \in \hL^0_{\hF^Y}(\sC_T(\sP_1))$, such that  $\sT(\hat\m)=\hat\m$.
Then, denoting
$(\hat X, \hat L):=(X^{\hat\m}, L^{\hat\m})$ to be the corresponding solution to (\ref{SDE2}) with $\m=\hat\m$, and  $\hat\hP=\hP^{\hat\m}$,  the Kallianpur-Strieble formula (\ref{muth}) implies that
$\hat\m_t=\m_t^{\hat X|Y}$, $t\in[0,T]$, under $\hat\hP$. In other words,  writing
 $$\hat B^{2}_t=Y_t-\int_0^t h(s, \hat X_{\cd\wedge s}, Y_{\cd\wedge s}, \hat\m_{\cd\wedge s})ds,
 $$
we see that $(\O^0, \cF^0, \hF^0, \hat\hP, ( B^1, \hat B^{2}), (\hat X,Y))$  is a weak solution of (\ref{SDE1}).

In order to make our scheme to work  we shall carry out the following  tasks in the next sections.

(i) Identify a subspace $\sE\subset \hL^0_{\hF^Y}(\sC_T(\sP_1))$, and show that the solution mapping $\sT$ is from $\sE$ to itself.

(ii) Show that we can at least find a sequence of $\cF^Y$-stopping times $\{\t_N\}$, such that $\sT^N:=\sT\big|_{[0,\t_N]}$ is a contraction, hence
has a fixed point $\hat \m^N$ on $[0,\t_N]$. We then argue that these $\hat\m^N$'s can be ``patched" together to become a fixed point $\hat\m\in\sE$.

(iii) Show that the law of the solution $(X, Y)$ is unique.
\begin{rem}
\label{remark01}
{\rm We should note that the ``localization" procedure is merely technical, in order to deal with the unboundedness caused by the
fraction in (\ref{muth}). In fact, such a technicality only occurs in the CMVSDEs when $h\neq 0$, and it is the fundamental difference between the SDE (\ref{SDE1}) and the CMVSDEs of the common noise type that we often see in the literature.
\qed}
\end{rem}

To simplify notations, in what follows  for each $\m\in\hL^0_{\hF^Y}(\sC_T(\sP_1))$ and the corresponding solution to (\ref{SDE2}) $(X^\m, L^\m)$, we denote $L^{-\m}=(L^\m)^{-1}$. We note that $L^{-\m}$ is the inverse Girsanov kernel of $L^\m$, and it is a $\hP^\m$-martingale, but not a $\hQ^0$-martingale. Last but not least, for any $\xi\in \hC_T$, we shall also use the notation $\xi^*_t=\sup_{s\in[0,t]}|\xi_s|$, and for any $p>0$, $\xi^{ *, p}_t=[\xi^*_t]^p$, $t\in[0,T]$.

\section{The Main Estimates}
\setcounter{equation}{0}

In this section we establish the main estimates that will be crucial for us to implement the solution scheme.
%
Before we start, we emphasize again that the conditional law $\m^{X|Y}$ in CMVSDE (\ref{SDE1}) is under the probability $\hP$ (under which $(B^1, B^2)$ is a Brownian motion), but our scheme is defined under the reference measure $\hQ^0$, connected to $\hP$ via a Girsanov kernel $L$, defined by SDE (\ref{SDE2}) or explicitly by (\ref{Dade}).

Now, for any $\mu,\mu'\in \hL^0_{\hF^Y}([0,T];\sP_1)$, denote
$\widetilde{\mu}:=\sT(\mu)$,\ $\widetilde{\mu}':=\sT(\mu')$. Let $(X^\m, L^\m)$, $(X^{\m'}, L^{\m'})$ be the solution of SDE (\ref{SDE2})
and define
\bea
\label{zeta}
\zeta_t(\mu,\mu'):=\hE^{\hQ^0}[(L^\mu)_t^{*,4}+(L^{-\mu})_t^{*,4}+(L^{\m'})_t^{*,4}+(L^{-\m'})_t^{*,4}\big|\mathcal{F}_t^Y], \q t\in[0,T].
\eea
Then it is not hard to check that $\zeta(\m, \m')$ is a continuous, increasing $\mathbb{F}^Y$-adapted process with $\zeta_0(\mu,\mu')=4$.
We have the following result.
\begin{prop}
\label{yl2}
Assume that Assumption \ref{assump1} is in force. Then, for all $0\leq s\leq t\leq T$, it holds $\hQ^0$-almost surely that
\bea
\label{4}
W_1(\widetilde{\mu}_s,\widetilde{\mu}'_t)\leq C\zeta_t(\mu,\mu')\{(\hE^{\hQ^0}[|X_s^\mu-X_t^{\mu'}|^2\big|\mathcal{F}_T^Y])^{\frac{1}{2}}+(\hE^{\hQ^0}[|L_s^\mu-L_t^{\mu'}|^2\big|\mathcal{F}_T^Y])^{\frac{1}{2}}\},
\eea
where $\zeta(\m,\m')$ is defined by (\ref{zeta}). Furthermore,
for each $p\ge 1$, there exists a $C_p>0$, such that $ \hE^{\hQ^0}[\zeta_T^p(\mu,\mu')]\leq C_p$, for all $\mu,\mu'\in \hL^0_{\hF^Y}([0,T];\sP_1)$.
\end{prop}

\noindent{\it Proof.}
First recall the {\it Kantorovich-Rubinstein formula} (see (\ref{W1})): For $0\leq s\leq t\leq T$,
$$ W_1(\widetilde{\mu}_s,\widetilde{\mu}'_t)=\sup\Big\{\Big|\int_{\mathbb{R}}\varphi d\widetilde{\mu}_s-\int_{\mathbb{R}}\varphi d\widetilde{\mu}'_t \Big|, ~\varphi\in \mbox{Lip}_1(\mathbb{R})\Big\}. $$
Since both $\widetilde{\mu}_s,\widetilde{\mu}'_t$ are probability measures, it suffices to consider only those test functions $\varphi\in \mbox{Lip}_1(\mathbb{R})$ with $\varphi(0)=0$ so that $|\varphi(z)|\leq|z|,\ z\in\mathbb{R}$. In other words, we
can write:
\bea
\label{KR}
W_1(\widetilde{\mu}_s,\widetilde{\mu}'_t)=\sup\Big\{\Big|\int_{\mathbb{R}}\varphi d\widetilde{\mu}_s-\int_{\mathbb{R}}\varphi d\widetilde{\mu}'_t\Big|,~ \varphi\in \mbox{Lip}_1(\mathbb{R}),\ \varphi(0)=0\Big\}.
\eea

Note that for any $\varphi\in Lip_1(\mathbb{R})$ with $\varphi(0)=0$, by definition of $\tilde\m$,  $\tilde\m'$ (see (\ref{muth})) we have
$$ \Big|\int_{\mathbb{R}}\varphi d\widetilde{\mu}_s-\int_{\mathbb{R}}\varphi d\widetilde{\mu}'_t\Big|=\big|\hE^\hP[\varphi(X_s^\mu)\big|\mathcal{F}_s^Y]-\hE^{\hP'}[\varphi(X_t^{\mu'})\big|\mathcal{F}_t^Y]\big|,$$
where $\hP=\hP^\mu$ and $\hP'=\hP^{\mu'}$. By the Kalliapur-Strieble formula (\ref{muth}) we have
$$
\hE^\hP[\varphi(X_s^\mu)\big|\mathcal{F}_s^Y]=\frac{\hE^{\hQ^0}[L_s^\mu \varphi(X_s^\mu)\big|\mathcal{F}_s^Y]}{\hE^{\hQ^0}
[L_s^\mu\big|\mathcal{F}_s^Y]}=\frac{\hE^{\hQ^0}[L_s^\mu \varphi(X_s^\mu)\big|\mathcal{F}_T^Y]}{\hE^{\hQ^0}[L_s^\mu\big|\mathcal{F}_T^Y]},~ \hQ^0\mbox{-a.s.},
$$
where the second equality follows from the fact that $\cF^Y_T=\cF^Y_s\vee\cF^Y_{s,T}$, and $\mathcal{F}_s=\mathcal{F}_s^{B^1,Y}$ and $\mathcal{F}_{s, T}^Y$ are  independent under $\hQ^0$. Similarly, we have
$$
\hE^{\hP'}[\varphi(X_t^{\mu'})\big|\mathcal{F}_t^Y]=\frac{\hE^{\hQ^0}[L_t^{\mu'} \varphi(X_t^{\mu'})\big|\mathcal{F}_t^Y]}{\hE^{\hQ^0}[L_t^{\mu'}\big|\mathcal{F}_t^Y]}=\frac{\hE^{\hQ^0}[L_t^{\mu'} \varphi(X_t^{\mu'})\big|\mathcal{F}_T^Y]}{\hE^{\hQ^0}[L_t^{\mu'}\big|\mathcal{F}_T^Y]}, ~\hQ^0\mbox{-a.s.}
$$
Hence, we deduce that, $\hQ^0$-almost surely,
\bea
\label{fest}
&&\Big|\int_{\mathbb{R}}\varphi d\widetilde{\mu}_s-\int_{\mathbb{R}}\varphi d\widetilde{\mu}'_t\Big|=\Big|\frac{\hE^{\hQ^0}[L_s^\mu \varphi(X_s^\mu)\big|\mathcal{F}_T^Y]}{\hE^{\hQ^0}[L_s^\mu\big|\mathcal{F}_T^Y]}-\frac{\hE^{\hQ^0}[L_t^{\mu'} \varphi(X_t^{\mu'})\big|\mathcal{F}_T^Y]}{\hE^{\hQ^0}[L_t^{\mu'}\big|\mathcal{F}_T^Y]}\Big|\nonumber\\
&\leq& \frac{1}{\hE^{\hQ^0}[L_s^\mu\big|\mathcal{F}_T^Y]}\hE^{\hQ^0}[|L_s^\mu\varphi(X_s^\mu)-L_t^{\mu'}\varphi(X_t^{\mu'})|\big|\mathcal{F}_T^Y]\\
&&\ \   +\frac{|\hE^{\hQ^0}[L_t^{\mu'}\varphi(X_t^{\mu'})\big|\mathcal{F}_T^Y]|}{\hE^{\hQ^0}[L_t^{\mu'}\big|\mathcal{F}_T^Y]\hE^{\hQ^0}[L_s^{\mu}\big|\mathcal{F}_T^Y]}\hE^{\hQ^0}[|L_s^\mu-L_t^{\mu'}|\big|\mathcal{F}_T^Y]=: I_{s,t}^1+I_{s,t}^2, \nonumber
\eea
where $I^i_{s,t}$, $i=1,2$, are defined in an obvious way.
Now by Jensen's inequality we have (recall the definition of $L^{-\m}$),
\bea
\label{Lmest}
(\hE^{\hQ^0}[L_s^{\mu}\big|\mathcal{F}_T^Y])^{-1}\leq \hE^{\hQ^0}[L_s^{-\mu}\big|\mathcal{F}_T^Y] \q\mbox{and}\q (\hE^{\hQ^0}[L_t^{\mu'}\big|\mathcal{F}_T^Y])^{-1}\leq \hE^{\hQ^0}[L_t^{-\mu'}\big|\mathcal{F}_T^Y],
\eea
and recalling the notation $\xi^*$ for $\xi\in\hC_T$, we have
\bea
\label{zeta1}
&&\Big|\frac{\hE^{\hQ^0}[L_t^{\mu'}\varphi(X_t^{\mu'})\big|\mathcal{F}_T^Y]}{\hE^{\hQ^0}[L_t^{\mu'}\big|\mathcal{F}_T^Y]\hE^{\hQ^0}[L_s^{\mu}\big|\mathcal{F}_T^Y]}\Big|\\
&\leq&\big(\hE^{\hQ^0}[(L^{\mu'})_t^{*,2}\big|\mathcal{F}_T^Y]\big)^{\frac{1}{2}}\big(\hE^{\hQ^0}[(X^{\mu'})_t^{*,2}\big|\mathcal{F}_T^Y]\big)^{\frac{1}{2}}\hE^{\hQ^0}[(L^{-\mu'})_t^*\big|\mathcal{F}_T^Y]\hE^{\hQ^0}[(L^{-\mu})_t^*\big|\mathcal{F}_T^Y]=:\zeta_t^1.
\nonumber
\eea
To analyze $\zeta^1$ we first recall that, under $\hQ^0$, $(B^1, Y)$ is a 2-dimensional Brownian motion. Therefore,
if we denote the conditional probability $\hQ^0[A|\cF^Y_T](\o^2)=\hQ^{\o^2}[A]$, $A\in \sB(\hC^2_T)$, then we can consider
the SDE for $X^\m$ in (\ref{SDE2}) as on the probability space $(\hC_T, \sB(\hC_T), \hQ^{\o^2})$ for $\hP^0$-a.e. $\o^2\in \hC_T$. Note that for fixed $\o^2$, the process
$$X^\m_t(\cd, \o^2)=x+\int_0^t\si(s, X^\m_s(\cd, \o^2), \o^2_{\cd\wedge s}, \m_{\cd\wedge s}(\o^2))dB^1_s, \q t\in[0,T],
$$
is an $\hQ^{\o^2}$-martingale, and as $\sigma$ is bounded, by the Burkholder-Davis-Gundy inequality we have
$$\hE^{\o^2}[(X^\m)^{*,2}_T]\le C \hE^{\o^2}[{\lan X^\m\ran}_T]=C\hE^{\o^2}\Big[\int_0^T\si^2(\cds)ds\Big]\le C, \q \hP^0\mbox{-a.e.
$\o^2\in \hC_T$,}
$$
where $\hE^{\o^2}[\, \cd\,]=\hE^{\hQ^{\o^2}}[\, \cd \, ]=\hE^{\hQ^0}[\, \cd\, |\cF^Y_T](\o^2)$, and $C>0$ is a generic
constant depending only on $T>0$ and the bounds of $\si$ and $h$, which is allowed to vary from line to line. Thus we have
$\hE^{\hQ^0}[(X^{\mu})_T^{*,2}\big|\mathcal{F}_T^Y]\leq C$, $\hQ^0$-a.s.

Now repeatedly applying H\"older's inequality and
the fact $abc\le a^3+b^3+c^3$, for $a,b,c\geq0$, we obtain from the definition of $\zeta_t^1$ in (4.6)
\bea
\label{zeta2}
\zeta_t^1&\le& C \hE^{\hQ^0}[(L^{\mu'})_t^{*,3}+(L^{-\mu'})_t^{*,3}+(L^{-\mu})_t^{*,3}\big|\mathcal{F}_t^Y]\\
&=& C \hE^{\hQ^0}[(L^{\mu'})_t^{*,3}+(L^{-\mu'})_t^{*,3}+(L^{-\mu})_t^{*,3}\big|\mathcal{F}_T^Y]=:
\zeta_t^2,\q t\in[0,T]. \nonumber
\eea
Now, for notational simplicity we denote $\D L_{s,t}^{\m,\m'}:=L^\m_s-L^{\m'}_t$, and $\D X_{s,t}^{\m,\m'}:=X^\m_s-X^{\m'}_t$.
Then, combining (\ref{Lmest})--(4.7), and recalling the definition of $I_{s,t}^2$ (see (\ref{fest})), we have
\bea
\label{Ist2}
I_{s,t}^2\leq\zeta_t^2 \hE^{\hQ^0}[|\D L_{s,t}^{\m,\m'}|\big|\mathcal{F}_T^Y],\q \hQ^0\text{-}a.s., ~0\leq s\leq t\leq T.
\eea
Similarly, we have the estimate for $I_{s,t}^1$ (noting that $(L^{\mu'})_t^*\geq L_0^{\mu'}=1$), for $0\leq s\leq t\leq T$,
\bea
\label{Ist1}
I_{s,t}^1\neg&\neg\leq\neg& \neg\hE^{\hQ^0}[(L^{-\mu})_t^*\big|\cF_T^Y]\big\{\big(\hE^{\hQ^0}[(X^{\mu})_t^{*,2}\big|\cF_T^Y]\big)^{\frac{1}{2}}
\big(\hE^{\hQ^0}[|\D L_{s,t}^{\m,\m'}|^2\big|\mathcal{F}_T^Y]\big)^{\frac{1}{2}}\nonumber\\
&&+\big(\hE^{\hQ^0}[(L^{\mu'})_t^{*,2}\big|\mathcal{F}_T^Y]\big)^{\frac{1}{2}}\big(\hE^{\hQ^0}[|\D X_{s,t}^{\m,\m'}|^2\big|
\mathcal{F}_T^Y]\big)^{\frac{1}{2}}\big\}\\
\neg&\neg\leq\neg& \neg C \hE^{\hQ^0}[(L^{-\mu})_t^{*,2}+(L^{\mu'})_t^{*,2}\big|\mathcal{F}_T^Y]\big\{(\hE^{\hQ^0}[|\D L_{s,t}^{\m,\m'}|^2\big|\cF_T^Y])^{\frac{1}{2}}+(\hE^{\hQ^0}[|\D X_{s,t}^{\m,\m'}|^2\big|\cF_T^Y])^{\frac{1}{2}}\big\}. \nonumber
\eea

Plugging (\ref{Ist2}) and (\ref{Ist1}) into (\ref{fest}) we have, for all $0\leq s\leq t\leq T$, $\hQ^0$-a.s.,
\bea
\label{5}
&&\Big|\int_{\mathbb{R}}\varphi d\widetilde{\mu}_s-\int_{\mathbb{R}}\varphi d\widetilde{\mu}'_t\Big|\\
&\leq& C \zeta_t(\mu,\mu')\Big((\hE^{\hQ^0}[|\D X_{s,t}^{\m,\m'}|^2\big|\mathcal{F}_T^Y])^{\frac{1}{2}}+
(\hE^{\hQ^0}[|\D L_{s,t}^{\m,\m'}|^2\big|\mathcal{F}_T^Y])^{\frac{1}{2}}\Big), \nonumber
\eea
where $C>0$ is a constant depending only on $T$ and the bounds of $\si$, $h$, and
\bea
\label{zeta2}
\zeta_t(\mu,\mu')&:=&\hE^{\hQ^0}[(L^{\mu})_t^{*,4}+(L^{-\mu})_t^{*,4}+(L^{\mu'})_t^{*,4}+(L^{-\mu'})_t^{*,4}\big|\cF_t^Y]\\
&=&\hE^{\hQ^0}[(L^{\mu})_t^{*,4}+(L^{-\mu})_t^{*,4}+(L^{\mu'})_t^{*,4}+(L^{-\mu'})_t^{*,4}\big|\mathcal{F}_T^Y],\ t\in[0,T]. \nonumber
\eea
From its definition we can easily see that $\zeta_t(\m,\m')$, $t\in[0,T]$, is an $\mathbb{F}^Y$-adapted,  increasing process
with $\zeta_0(\mu,\mu')=4$. Moreover, by the last expression of (\ref{zeta2}) we see that it is $L^2 (\hQ^0)$-continuous. Thus, the continuity of $t\rightarrow \zeta_t(\m,\m')$ follows. Finally, for each $p\ge 1$, there exists some constant $C_p>0$,
depending only $p$ and the bounds of coefficients, such that
$$ \hE^{\hQ^0}[\zeta_T^p(\mu,\mu')]\leq C_p,    \ \mbox{ for all } \mu,\ \mu'\in \hL_{\hF^Y}^{0}( [0,T]; \sP_1),\ p\geq1. $$
This proves the proposition.
\qed

We now consider the following subspace of $\hL^0_{\hF^Y}(\sC_T(\sP_1))$ (see (\ref{Sp}) for definition):
\bea
\label{sE*}
\sE:= \hS^{\infty-}_{\hF^Y}(\sP_1).
\eea
We shall argue that
%
the conclusion of Proposition \ref{yl2} is strong enough to imply the following important property of the solution mapping $\sT$.
\begin{cor}
\label{yl1}
Assume Assumption \ref{assump1}. Then $\sT(\sE)\subseteq \sE$.
\end{cor}
\noindent {\it Proof.}
For any $\m\in\sE$ we put $\widetilde{\mu}=\sT(\mu)$. Setting $\mu'=\mu$ in Proposition \ref{yl2}, we deduce from
 (\ref{4}) that
\bea
\label{must}
\hE^{\hQ^0}[W_1(\widetilde{\mu}_s,\widetilde{\mu}_t)^4]&\leq& C\big(\hE^{\hQ^0}[\zeta_T(\mu,\mu)^8]\big)^\frac{1}{2}
\Big\{\big(\hE^{\hQ^0}[|X^\mu_s-X^\mu_t|^8]\big)^\frac12+\big(\hE^{\hQ^0}[|L^\mu_s-L^\mu_t|^8]\big)^\frac{1}{2}\Big\} \nonumber\\
&\leq& C\Big\{\big(\hE^{\hQ^0}[|X^\mu_s-X^\mu_t|^8]\big)^\frac12+\big(\hE^{\hQ^0}[|L^\mu_s-L^\mu_t|^8]\big)^\frac{1}{2}\Big\}.
\eea
Here and in what follows we shall denote $C>0$ to be a generic constant depending only on $T$ and the bounds of $h$, which varies from line to line. Since $\sigma$ and $h$ are bounded, it is clear that $\hE^{\hQ^0}[(L^\mu_T)^p]\leq C_p$, for all $\mu\in \sC_T(\sP_1)$, and it follows by standard estimates that
$$\hE^{\hQ^0}[|X^\mu_s-X^\mu_t|^8+|L^\mu_s-L^\mu_t|^8]\leq C|s-t|^4,\quad 0\leq s\leq t\leq T.
$$
Hence, by (\ref{must}) we have $\hE^{\hQ^0}[W_1(\widetilde{\mu}_s, \widetilde{\mu}_t)^4]\leq C|s-t|^2$, $s, t\in[0,T]$. Thus,
by Kolmogorov's continuity criterion, it follows that $\widetilde{\mu} = (\widetilde{\mu}_t)_{t\in[0,T]}$ admits a continuous modification,
which we shall use from now on. In other words, we have proved that $\sT(\m)=\widetilde\m$ is $\sC_T(\sP_1)$-valued.

It remains to check that $\sT(\m)\in\sE$. To see this we fix $p \geq 1$, and note that for any $\m'\in \sE$ we always have
$\sT(\m')_0=\hP'\circ (X^{\m'}_0)^{-1}=\delta_{\{x\}}$. Applying  (\ref{4}) again we see that
\beaa
&&\hE^{\hQ^0}\Big[\sup_{t\in[0,T]}W_1(\widetilde{\mu}_t,\delta_{x_0})^p\Big] = \hE^{\hQ^0}\Big[\sup_{t\in[0,T]}W_1(\widetilde{\mu}_t,\widetilde{\mu}'_0)^p\Big]\\
&\leq& C\hE^{\hQ^0}\Big[\zeta_T(\mu,\mu')^p\cdot\sup_{t\in[0,T]}\Big(\big(\hE^{\hQ^0}[|X^\mu_t-x|^2|\mathcal{F}^Y_T]\big)^\frac{p}{2}+(\hE^{\hQ^0}[|L^\mu_t-1|^2|\mathcal{F}^Y_T]\big)^\frac{p}{2}\Big)\Big]\\
&\leq &C\Big(\hE^{\hQ^0}\big[\zeta_T(\mu,\mu')^{2p}\big]\Big)^\frac{1}{2}\Big(1+\hE^{\hQ^0}[(X^\mu)^{*, 2p}_T ]+\hE^{\hQ^0}[(L^\mu)^{*,2p}_T\big]\Big)^\frac{1}{2} < +\infty.
\eeaa
Since $\widetilde{\mu}$ is obviously $\mathbb{F}^Y$-adapted, by definition (\ref{Sp}) we then have $\widetilde\m\in \hS^p_{\hF^Y}(\sP_1)$. Note now that the above argument holds for all $p\geq1$,
we conclude that $\sT(\m)=\widetilde\m\in \sE$. The proof is now complete.
\qed

\begin{rem}
\label{remark3}
{\rm
As we pointed out before, Proposition \ref{yl2} actually shows that $\sT(\m)\in \hS^{\infty-}_{\hF^Y}(\sP_1)$, for any
$\m\in \hL^0_{\hF^Y}([0,T];\sP_1)$. This is due largely to the fact that we have assumed that all coefficients $\si$ and $h$ are bounded. In general, we should have, for any $p\ge 1$, the solution mapping $\sT:\hS^p_{\hF^Y}(\sP_1)\mapsto  \sE=
\hS^{\infty-}_{\hF^Y}(\sP_1)\subseteq \hS^p_{\hF^Y}(\sP_1)$. The case when $p=2$ is frequently used.
\qed}
\end{rem}

\section{Existence of weak solution}
\setcounter{equation}{0}

We are now ready to prove the existence of the weak solution to SDE (\ref{SDE1}).
To begin with, we note that Proposition \ref{yl2} only shows that (assuming, for example, $s=t$), the (Wasserstein) distance between
$\tilde\m_t=\sT(\m)_t$ and $\tilde\m'_t=\sT(\m')_t$ can be controlled by the distances of the corresponding solutions $(X^\mu, L^\mu)$
and $(X^{\mu'}, L^{\mu'})$ at each fixed $t\in[0,T]$. But in order to look for a fixed point in the space $\sE$, we need to strengthen the estimate in terms of the distance in $\hL^p_{\hF^Y}(\sC(\sP_1))$. In light of Remark \ref{remark3}, we shall only consider the case
$p=2$.


We begin by a brief analysis. Let $\mu, \mu'\in \hS^2_{\hF^Y}(\sP_1)$. For notational simplicity we denote the corresponding triplets
$(X,L,\hP):=(X^\mu,L^\mu,\hP^\mu)$ and $(X',L',\hP'):=(X^{\mu'},L^{\mu'},\hP^{\mu'})$, respectively, and put $\widetilde{\mu}:=\sT(\mu)$ and $\widetilde{\mu}':=\sT(\mu')$ as before. We also set $\D X :=X-X'$, $\D L :=L-L'$, and
$$ \d \f(t, x, x', y, \m,\m'):=\f(t,x,y, \m)-\f(t, x',y, \m'), \q \f=\si,\ h.
$$
Our goal is to use estimate (\ref{4}) in Proposition \ref{yl2} to obtain the desired contraction estimate: For some constant $C\in (0, 1)$,
\bea
\label{contra}
\hE^{\hQ^0}\Big[\sup_{0\le t\le T}W_1(\tilde \m,\tilde \m')^2\Big]\le C\hE^{\hQ^0}\Big[\sup_{0\le t\le T}W_1(\m,\m')^2\Big].
\eea

To begin with, we note that (\ref{4}) only gives us
\bea
\label{8}
\sup_{s\leq t} W_1(\widetilde{\mu}_s,\widetilde{\mu}'_s)^2\leq C\zeta_t(\mu,\mu')^2 \hE^{\hQ^0}\Big[\sup_{s\leq t}|\D X_s|^2 +
\sup_{s\leq t}|\D L_s|^2 \big|\mathcal{F}_T^Y\Big],\quad t\in[0,T].
\eea

But on the other hand, since $X$ and $X'$ satisfy (\ref{SDE2}), following the  standard arguments using the Burkholder-Davis-Gundy inequality and Assumption \ref{assump1}, one can easily check that,
for $t\in[0,T]$,
\beaa
&&\hE^{\hQ^0}[\sup_{s\in[0,t]}|\D X_s|^4\big|\mathcal{F}_T^Y]
\leq \ C\hE^{\hQ^0}\Big[\Big(\int^t_0\big|\d\sigma(s,X_{\cd\wedge s},X'_{\cdot\wedge s},Y_{\cdot\wedge s},\mu_{\cdot\wedge s},
\mu'_{\cdot\wedge s})\big|^2ds\Big)^2\Big|\mathcal{F}_T^Y\Big]\\
&\leq& C\hE^{\hQ^0}\Big[\int^t_0 \sup_{r\leq s}|\D X_r|^4 ds\Big|\mathcal{F}_T^Y\Big]
+C\hE^{\hQ^0}\Big[\Big(\int^t_0 \sup_{r\leq s} W_1(\mu_r,\mu'_r)^2 ds\Big)^2\Big|\mathcal{F}_T^Y\Big]\\
&=& C\int^t_0 \hE^{\hQ^0}[\sup_{r\leq s}|\D X_r|^4 |\mathcal{F}_T^Y]ds
+C\Big(\int^t_0 \sup_{r\leq s} W_1(\mu_r,\mu'_r)^2 ds\Big)^2.
\eeaa
Observe that in the last equality above we used  the fact that $\mu$ and $\mu'$ are $\hF^Y$-adapted.
Now applying Gronwall's inequality  we obtain that
\bea
\label{6}
\Big(\hE^{\hQ^0}\Big[\sup_{s\in[0,t]}|\D X_s|^4 \big|\mathcal{F}_T^Y \Big]\Big)^\frac{1}{2}
\leq C\int^t_0 \sup_{r\leq s} W_1(\mu_r,\mu'_r)^2 ds,\  t\in[0,T],\q \hQ\mbox{-}a.s.
\eea

Similarly, since $h$ is bounded, we also obtain from (\ref{SDE2}) that, for $t\in[0,T]$,
\beaa
 \hE^{\hQ^0}\Big[\sup_{s\leq t}|\D L_s|^2 \Big]&\leq& C\Big\{\hE^{\hQ^0}\Big[\int^t_0 |\D L_s|^2 ds\Big]\\
&&\ \ \ \ + \hE^{\hQ^0}\Big[\int^t_0 L^2_s\big|\d h(s,X_{\cdot\wedge s},X'_{\cdot\wedge s},Y_{\cdot\wedge s},\mu_{\cdot\wedge s},\mu'_{\cdot\wedge s})\big|^2 ds\Big]\Big\},
\eeaa
and again applying Gronwall's inequality we get
\beaa
&&\hE^{\hQ^0}\Big[\sup_{s\in[0,t]}|\D L_s|^2\Big]
\leq\ C\hE^{\hQ^0}\Big[\int^t_0 L^2_s\Big(\sup_{r\leq s}|\D X_r|^2
+ \sup_{r\in[0,s]} W_1(\mu_r,\mu'_r)^2\Big) ds\Big]\\
&\leq& C\hE^{\hQ^0}\Big[\int^t_0 \Big(\hE^{\hQ^0}\big[L^4_s\big|\mathcal{F}_T^Y\big]\Big)^\frac{1}{2}\Big[\Big(\hE^{\hQ^0}\Big[\sup_{r\leq s}|\D X_r|^4 \big|\mathcal{F}_T^Y\Big]\Big)^\frac{1}{2}+ \sup_{r\leq s} W_1(\mu_r,\mu'_r)^2 \Big]ds\Big].
\eeaa
Now by (\ref{6}) we conclude from the above that, for $t\in[0,T]$, $\hQ^0$-a.s.,
\bea
\label{7}
\hE^{\hQ^0}\Big[\sup_{s\in[0,t]}|\D L_s|^2\Big]
\leq\ C\hE^{\hQ^0}\Big[\int^t_0 \big(\hE^{\hQ^0}\big[L^4_s\big|\mathcal{F}_T^Y\big]\big)^\frac{1}{2}\cdot \sup_{r\leq s} W_1(\mu_r,\mu'_r)^2 ds\Big].
\eea

Moreover, noting that $\hE^{\hQ^0}\big[L^4_s\big|\mathcal{F}_T^Y\big]\le \zeta_t(\m,\m')$, we see from (\ref{8}), (\ref{6}), and (\ref{7})   that we would easily have the desired estimate (\ref{contra}) so the Contraction Mapping Theorem can be applied (at least in the case when the time duration is small) if we could find a bound for $\zeta_t(\mu,\mu')$
that is independent of $\mu$, $\mu'$. But this is in general difficult, since each $L^\m$ is the solution to a linear SDE driven by
the $\hQ^0$-Brownian motion $Y$, thus under the conditional expectation $\hE^{\hQ^0}[\cd|\cF^Y_t]$, this essentially amounts to asking a pathwise uniform bound for a family of martingales, which is generally impossible. We shall therefore impose the following extra structural
assumption on the coefficient $h$ in SDE (\ref{SDE1}).
\begin{assum}
\label{assump2}
The function $h$ in (\ref{SDE1}) is of the form:
\bea
\label{h}
 h(t,x,y_{\cdot\wedge t}) = \sum^N_{i=1} f_i(t, x)g_i(t,y_{\cdot\wedge t}),
 \eea
 where $f_i\in \hC^{1, 2}_b([0,T]\times \hR)$, $1\le i\le N$, and $g_i$'s are bounded and measurable.
\qed
 \end{assum}

We remark that the Assumption \ref{assump2}  trivially contains all the traditional nonlinear filtering problems, in which $h=h(t,x)$.
In what follows, without loss of generality we shall assume
$N = 1$, and $f:=f_1,\ g:=g_1$. We have the following crucial result regarding the process $\zeta(\m,\m')$ defined by (\ref{zeta}), for any $\mu, \mu' \in \hS_{\mathbb{F}^Y}^2(\sP_1)$.
\begin{prop}
\label{zetaest}
Assume the Asumptions \ref{assump1} and \ref{assump2} are in force. Then there exists a continuous,  increasing, $\hF^Y$-adapted process $A=\{A_t\}_{t\in[0,T]}$, with
$A_0 >0$,
such that for any $\mu, \mu' \in S_{\mathbb{F}^Y}^2(\sP_1)$, it holds that  $\zeta_t(\mu, \mu')\leq A_t$, $t\in[0,T]$, $\hQ^0$-a.s.
\end{prop}

\noindent{\it Proof.}
For any $\m\in \hS_{\hF^Y}^2(\sP_1)$, let $X=X^\m$ be the solution to (\ref{SDE2}).  Since $f\in \hC^{1,2}([0,T]\times\hR)$,  thanks to Assumption \ref{assump2}, we apply It\^o's formula to get:
\bea
\label{Ito}
df(t, X_t)= \pa_tf(t, X_t)dt+\pa_xf(t,X_t)dX_t+\frac{1}{2}\pa^2_{xx}f(t,X_t)d\langle X\rangle_t, \q t\in[0,T].
\eea

Now let us consider the following two processes:
\bea
\label{ZM}
\left\{\ba{lll}
\dis Z_t=\int_0^tg(s,Y_{\cd\wedge s})dY_s;\ms\\
\dis M_t=M^\m_t:=\int_0^tZ_s\pa_xf(s, X_s)\sigma(s, X_{\cdot\wedge s}, Y_{\cdot\wedge s}, \mu_{\cdot\wedge s})dB_s^1,
\ea\right. \q t\in[0,T].
\eea
%
Recalling that $B^1$ and $Y$ are independent under $\hQ^0$, we have
$d{\lan M, Z\ran}_t=d\lan X, Z\ran_t\equiv 0$.
%
Thus (\ref{Ito}) and  integrating by parts yield, for $t\in[0,T]$,
\bea
\label{hest}
&&\int_0^th(s, X_s, Y_{\cdot\wedge s})dY_s=\int_0^tf(s, X_s)g(s, Y_{\cdot\wedge s})dY_s=\int_0^tf(s,X_s)dZ_s\nonumber\\
&=&f(t, X_t)Z_t-\int_0^tZ_s\big[\pa_tf(s, X_s)ds+\pa_xf(t,X_s)dX_s+\frac{1}{2}\pa^2_{xx}f(s,X_s)d\langle X\rangle_s\big]\\
&=&f(t,X_t)Z_t-M_t-\int_0^tZ_s\big[\pa_tf(s,X_s)+\frac12\pa_{xx}^2f(s,X_s)|\sigma(s, X_{\cdot\wedge s}, Y_{\cdot\wedge s},
\mu_{\cdot\wedge s})|^2\big]ds.
 \nonumber
\eea
Since $\si$ is bounded and $f\in \hC^{1,2}_b$  we easily deduce that
\bea
\label{Zstar}
\left\{\ba{lll}
\dis \Big|f(t,X_t)Z_t-\int_0^tZ_s\big[\pa_tf(s,X_s)+\frac12\pa_{xx}^2f(s,X_s)|\sigma(s, X_{\cdot\wedge s}, Y_{\cdot\wedge s},
\mu_{\cdot\wedge s})|^2\big]ds\Big|
\le CZ_t^*, \ms \\
\dis {\lan M\ran}_t=\int_0^t |Z_s\pa_x f(s,X_s)\si(s, X_{\cdot\wedge s}, Y_{\cdot\wedge s}, \m_{\cd\wedge s})|^2ds \le CZ^{*,2}_t,  \q t\in[0,T].
\ea\right.
\eea
Here $C>0$ is a generic constant depending only on the bounds of $f$ and $g$. Since $Z^*$ is  $\hF^Y$-adapted, a direct computation
using (\ref{hest}) and (\ref{Zstar}) shows  that, for all $p>0$, $t\in[0,T]$,
\bea
\label{phest}
&&\hE^{\hQ^0}\Big[\sup_{s\leq t}\Big(\exp\Big\{p\int_0^sh(r, X_r, Y_{\cdot\wedge r})dY_r\Big\}\Big)\Big|\cF_T^Y\Big]
\nonumber\\
\neg\neg&\neg\neg=\neg\neg&\neg\hE^{\hQ^0}\Big[\sup_{s\leq t}\neg\Big(\neg\exp\Big\{pf(s, X_s)Z_s-pM_s-p\int_0^sZ_r[\pa_tf(r,X_r)\nonumber\\
&&\qq\qq\qq\ \ +\frac{1}{2}\pa_{xx}^2f(r,X_r)|\si(r,X_{\cdot\wedge r}, Y_{\cdot\wedge r}, \m_{\cd\wedge r})|^2]dr\Big\}\Big)\Big|\mathcal{F}_T^Y\Big]\nonumber\\
\neg&\neg\leq\neg&\neg\neg C_p\Big(\hE^{\hQ^0}\big[\sup_{s\leq t}\big(\exp\{-pM_s-p^2\langle M\rangle_s\}\big)^2|\cF_T^Y\big]\Big)^{\frac{1}{2}}e^{C_pZ_t^{*,2}}\\
\neg&\neg\leq\neg&\neg C_p\Big(\hE^{\hQ^0}[\exp\{-2pM_t-2p^2\langle M\rangle_t\}|\mathcal{F}_T^Y]\Big)^{\frac{1}{2}}e^{C_pZ_t^{*,2}}
= C_p\big(\hE^{\hQ^0}[\cE_t|\cF_T^Y]\big)^{\frac{1}{2}}e^{C_pZ_t^{*,2}},\nonumber
\eea
where $C_p>0$ is some generic constant that may depends on $p$, and is allowed to vary from line to line,
and $\cE_t:=\exp\{-2pM_t-\frac12\langle 2pM\rangle_t\}$ is the Dol\'eans-Dade stochastic exponential
of the process $2pM$. That is, $\cE$ solves the linear SDE:
\bea
\label{cE}
\cE_t&=&1-\int_0^t \cE_sd(2pM_s)\nonumber\\
&=&1-2p\int_0^t \cE_sZ_s\pa_xf(s,X_s)\si(s,X_{\cdot\wedge s}, Y_{\cdot\wedge s}, \m_{\cd\wedge s})dB^1_s, \q t\in[0,T].
\eea
Now
consider the regular conditional probability
$\hP^{\o^2}_{T}(\cd):=\hQ^0[\,\cd\,|\cF^Y_T](\o^2)$, for $\hP_0$-a.e. $\o^2\in\hC_T$. For an $\hF$-adapted process $\xi$ we denote
$\xi^{\o^2} (\o^1)=\xi (\o^1, \o^2)$, $(\o^1, \o^2)\in\hC_T^2$. Then, since $\m$ is $\hF^Y$-adapted, (\ref{cE})  means that for
  $\cE^{\o^2} $,
for $\hP_0$-a.e. $\o^2\in\hC_T$, it holds
$\hP^{\o^2}_{T}$ almost surely:
\beaa
\label{cE1}
\cE_t^{\o^2}=
1-2p\int_0^t \cE^{\o^2}_sZ_s^{\o^2}\pa_xf(s,X^{\o^2}_s)\si(s, X^{\o^2}_{\cdot\wedge s},  \o^2_{\cd\wedge s},  \m^{\o^2}_{\cdot\wedge s})dB^1_s, \q t\in[0,T].
\eeaa
 That is, $\cE^{\o^2}$ is an exponential martingale under $\hP^{\o^2}_{T}$, and thus
 $\hE^{\hQ^0}[\cE_t|\cF^Y_T](\o^2)=1$, for $\hP_0$-a.e. $\o^2\in\hC_T$.
Consequently, since $h$ is bounded, for $\mu\in \hS^2_{\hF^Y}(\sP_1)$ and $p\ge 1$ we see from the definition of $L^\m$ (\ref{Dade}) and  (\ref{phest}) that
\beaa
\hE^{\hQ^0}[(L^\m)_t^{*, p}|\cF_T^Y] \leq  C_p\hE^{\hQ^0}\Big[\sup_{0\le s\leq t}\Big(\exp\Big\{p\int_0^sh(r, X_r, Y_{\cdot\wedge r})dY_r\Big\}\Big)\Big|\cF_T^Y\Big]
\le C_p e^{C_pZ_t^{*,2}}.
\eeaa
But this particularly implies that there exists a constant $C>0$, such that for any $\mu,\mu'\in \hS_{\mathbb{F}^Y}^2(\sP_1)$,
it holds that
\bea
\label{A}
\zeta_t(\mu, \mu')\leq
C\exp{\{C Z^{*,2}_t\}}=:A_t,  \q t\in[0,T],
\eea
where $Z$ is defined by (\ref{ZM}). Clearly, the process $A$ is   continuous, $\hF^Y$-adapted, increasing, and is  independent of
the choice of $\m$, proving the proposition.
\qed

We now give the main result of this section.
\begin{thm}
\label{fixedpt}
Assume the Assumptions \ref{assump1} and  \ref{assump2}. Then, the solution mapping $\sT(\cd)$ defined by (\ref{muth}) has a unique fixed point in $ \hS_{\mathbb{F}^Y}^2(\sP_1)$.
\end{thm}

\noindent{\it Proof.}  First consider the process $A$ in Proposition \ref{zetaest}. For $N\geq 1$, define the $\mathbb{F}^Y$-stopping time $\tau_N:=\inf\{t\geq0:A_t>N\}\wedge T$. Then, $\hQ^0\{\t_N\nearrow T\}=1$. Moreover,
let us now define, for $N\in \hN$ and $p\ge 1$, $ \hS_{\hF^Y}^{p,N}(\sP_1):=\{\mu_{\cdot\wedge\tau_N}:\mu\in \hS_{\hF^Y}^{p}(\sP_1)\}$, and
\bea
\label{13}
\sT_N(\mu)_t:=\sT(\mu)_{t\wedge\tau_N},  \q t\in[0,T], \q \mu\in \hS_{\hF^Y}^{p,N}(\sP_1).
\eea
Then, applying Proposition \ref{yl2} and Corollary \ref{yl1} we conclude that $\sT_N$ is a mapping from $\hS_{\hF^Y}^{2, N}(\sP_1)$ to itself.

We first show that each $\sT_N$, $N\in \hN$, has a fixed point.
To this end,  let $\m\in \hS_{\hF^Y}^{2, N}(\sP_1)$ and let $(X^\m, L^\m)$ be the corresponding solution of (\ref{SDE2}). Consider the
function $(t, x, \o)\mapsto h(t, x, Y_{\cd\wedge t}(\o))\1_{[0,\t_N]}(t,\o)$. Since $\tau_N$ is an $\mathbb{F}^Y$-stopping time, and $Y$ is a canonical Brownian motion ($Y_t(\o)=\o^2_t$) under $\hQ^0$, there is some  bounded and measurable functional $h^N:[0,T]\times\hR\times \hC_T\rightarrow \mathbb{R}$, such that

(i) for each $x\in \hR$, and $(t,\o)\in[0,T]\times\O^0$,
$ h^N(t, x, \o)=h^N(t, x, \o^2_{\cd\wedge t})$. In other words, the mapping $(t, \o)\mapsto h^N(t, x, \o)$ is $\hF^Y$-progressively measurable; and

(ii) it holds that
\bea
\label{hN}
h^N(t, X^\m_t, Y_{\cdot\wedge t})=h(t, X^\m_t, Y_{\cdot\wedge t})\1_{[0,\tau_N]}(t,\o),\ \  t\in [0,T].
\eea
Using the function $h^N$ we can solve the SDE:
\bea
\label{LmN}
L^{\m,N}_t=1+\int_0^th^N(s, X^\m_s, Y_{\cdot\wedge s})L^{\m,N}_sdY_s, \ \  t\in[0,T].
\eea
Then, by uniqueness, it is easy to check that $L^{\m, N}\equiv L^\m_{\cd\wedge\tau_N}$, where $L^{\m}$ solves (\ref{SDE2}).

Now for $\m, \m'\in \hS_{\hF^Y}^{2, N}(\sP_1)$, let $(X^\m, L^\m)$, $(X^{\m'}, L^{\m'})$ be the corresponding solutions to (\ref{SDE2}),
respectively. We shall denote $X:=X^\m$, $X':=X^{\m'}$,  and $L^N:=L^{\m,N}$, $L'^{N}:=L^{\m',N}$ for simplicity.
Since $ h^N$ is uniformly Lipschitz continuous in $x$ with the same Lipschitz constant as $h$ given in Assumption \ref{assump1}, we deduce from (\ref{6})-(\ref{7}) that, for $t\in[0,T]$,
\bea
\label{1.6'}
\hE^{\hQ^0}\Big[\sup_{s\in[0,t]}|X_s-X'_s|^4\Big|\mathcal{F}_T^Y\Big]&\leq& C\Big[\int_0^t\sup_{r\leq s}W_1(\mu_r, \mu'_r)^2ds\Big]^{2},
\\
\label{1.7'}
\hE^{\hQ^0}\Big[\sup_{s\in[0,t]}|L_s^N-L'^N_s|^2\Big]&\leq& C\hE^{\hQ^0}\Big[\int_0^t(\hE^{\hQ^0}[L_{s\wedge\tau_N}^4|\mathcal{F}_T^Y])^{\frac{1}{2}}\cd\sup_{r\leq {s\wedge \tau_N}}W_1(\mu_r, \mu'_r)^2ds\Big] \nonumber\\
&\leq& C\sqrt N\int_0^t\hE^{\hQ^0}\Big[\sup_{r\leq s\wedge{\tau_N}}W_1(\mu_r, \mu'_r)^2\Big]ds.
\eea
 Here in the above we used the facts that $\hE^{\hQ^0}[L_{s\wedge \tau_N}^4|\mathcal{F}_T^Y]\leq
 A_{t\wedge \tau_N}\leq N$ by  definition of $\t_N$; and that  $\tilde \m, {\tilde \m}'$ in (\ref{8}) should be replaced by
$\widetilde{\mu}^N_{t}(\cd):=\hP^N\{{{X_t}\in \cd}|\mathcal{F}_t^Y\}$ and $\widetilde{\mu}'^N_{t }(\cd):=\hP'^N\{X_t\in\cd|\mathcal{F}_t^Y\}$,  $t\in[0,T]$, where $d\hP^N:=L_T^Nd\hQ^0$, and $d\hP'^N:=L'^N_Td\hQ^0$, respectively.

On the other hand, recall again that $\tau_N$ is an $\mathbb{F}^Y$-stopping time, and for $B\in\mathcal{B}(\mathbb{R})$ we have
\beaa
\widetilde{\mu}_t^N(B)&=&\hP^N\{{{X_t}\in B}|\mathcal{F}_t^Y\}=\frac{\hE^{\hQ^0}[L^N_t\mathbf{1}_{\{X_t\in B\}}|\mathcal{F}_t^Y]}{\hE^{\hQ^0}[L_t^N|\mathcal{F}_t^Y]}=\frac{\hE^{\hQ^0}[L_{t\wedge\tau_N}\mathbf{1}_{\{X_t\in B\}}|\mathcal{F}_T^Y]}{\hE^{\hQ^0}[L_{t\wedge\tau_N}|\mathcal{F}_T^Y]}\\
&=&\frac{\hE^{\hQ^0}[L_t\mathbf{1}_{\{X_t\in B\}}|\mathcal{F}_T^Y]}{\hE^{\hQ^0}[L_t|\mathcal{F}_T^Y]}= \widetilde{\mu}_t(B), \q t\leq\tau_N,
\q \hQ^0\mbox{-}a.s.
\eeaa
In other words, we have $\tilde \m^N_{\cd\wedge \t_N}=\tilde \m_{\cd\wedge \t_N}=\sT_N(\m)$,  by definition (\ref{13}).  Similarly, we have
$\widetilde{\mu}'^N_{\cd\wedge\t_N}=\sT_N({\mu}')$.  Furthermore, if we modify (\ref{8}) as, for $t\in[0,T]$,
\bea
\label{9}
\sup_{s\leq t\wedge\tau_N}W_1(\widetilde{\mu}_s^N, \widetilde{\mu}'^N_s)^2\leq N^2\hE^{\hQ^0}\Big[\sup\limits_{s\leq t\wedge\tau_N}|X_s-X'_s|^2+\sup_{s\leq t\wedge\tau_N}|L_s-L'_s|^2\Big|\mathcal{F}_T^Y\Big],
\eea
then combining with (\ref{1.6'}) and (\ref{1.7'}) we derive,  for all $ \mu,\mu'\in S_{\mathbb{F}^Y}^{2,N}(\mathcal{P}_1)$, that
\bea
\label{14}
\hE^{\hQ^0}\Big[\sup_{s\leq t}W_1(\sT_N({\mu})_s, \sT_N({\mu}')_s)^2\Big]
&\leq& C_N\int_0^t\hE^{\hQ^0}\Big[\sup\limits_{r\leq s\wedge\tau_N}W_1(\mu_r, \mu'_r)^2\Big]ds\\
%
&\leq& C_N\int_0^t\hE^{\hQ^0}\Big[\sup_{r\leq s}W_1(\mu_r, \mu'_r)^2\Big]ds,\q t\in[0,T]. \nonumber
\eea
Consequently, iterating (\ref{14}) $k$ times we have, for $\mu,\ \mu'\in S_{\mathbb{F}^Y}^{2,N}(\mathcal{P}_1)$,
\beaa
\hE^{\hQ^0}\Big[\sup_{s\leq t}W_1\big(\sT_N^k(\mu)_s, \sT_N^k(\mu')_s\big)^2\Big]&\leq& C_N^k\int_0^t\neg\int_0^{t_1}\neg\neg\cdots\neg
\int_0^{t_{k-1}}\hE^{\hQ^0}\Big[\sup_{r\leq t_k}W_1(\mu_r,\mu'_r)^2\Big]dt_k\cdots dt_1\\
&\le&\frac{(C_Nt)^k}{k!}\hE^{\hQ^0}\Big[\sup_{r\leq t}W_1(\mu_r, \mu'_r)^2\Big], \q t\in[0,T].
\eeaa
Choosing $k$ large enough so that $(C_NT)^k/k!\leq\frac{1}{2}$, we see that $\sT_N^k:\ S_{\mathbb{F}^Y}^{2,N}(\sP_1)\rightarrow S_{\mathbb{F}^Y}^{2,N}(\sP_1)$ is a contraction, thus there is a unique $\bar\mu^N\in S_{\hF^Y}^{2,N}(\sP_1)$ such that
$\sT_N^k(\bar\mu^N)=\bar\mu ^N$, which implies that
$\sT_N(\bar\mu^N)=\sT_N(\sT_N^k(\bar\mu^N))=\sT_N^k(\sT_N(\bar\mu^N))$. Then by uniqueness we obtain $\sT_N(\bar\mu^N)=\bar\mu ^N$.

To construct the desired fixed point for $\sT$ on $[0,T]$ we shall argue that there is a standard extension, $\bar\m$, of the family $\{\bar\mu^N\}_{N\ge 1}$: $\bar\m_t=\bar\m^N_t$, whenever $t\in [0,\t_N]$, as $\hQ^0\{\t_N\nearrow T\}=1$. For this, we first claim that, given $ \mu\in S_{\hF^Y}^2(\sP_1)$ and any $\hF^Y$-stopping time $\tau\le T$, one has
\bea
\label{15}
\sT( \mu)_{t\wedge \tau}=\sT( \mu_{\cd\wedge \tau})_t, \qq t\in[0,T].
\eea
Indeed, for any bounded measurable function $\f$ and $t\in[0,T]$, we have
\beaa
&&\int_{\mathbb{R}}\varphi(x) \sT(\mu)_{t\wedge\tau}(dx)=\hE^{\hP^\mu}[\varphi(X_s^\mu)|\mathcal{F}_s^Y]\big|_{s=t\wedge\tau}={\frac{\hE^{\hQ^0}[\varphi(X_s^{\mu})L_s^\mu|\mathcal{F}_T^Y]}{\hE^{\hQ^0}[L_s^\mu|\mathcal{F}_T^Y]}}\bigg|_{s=t\wedge \tau}\\
&=&\frac{\hE^{\hQ^0}[\varphi(X_{t\wedge\tau}^{\mu})L_{t\wedge\tau}^{\mu}|\mathcal{F}_T^Y]}{\hE^{\hQ^0}[L_{t\wedge\tau}^{\mu}|\mathcal{F}_T^Y]}=\hE^{\hP^\mu}[\varphi(X^\m_{t\wedge\tau})|\mathcal{F}_{t\wedge\tau}^Y]=\int_{\mathbb{R}}\varphi(x) \sT({\mu}_{\cdot\wedge\tau})_t(dx).
\eeaa
This proves (\ref{15}). Now, using ({\ref{15}}) and the definition of $\sT_N$ we can further deduce that
$$
\sT_N(\bar\mu_{\cdot\wedge\tau_N}^{N+1})=\sT{(\bar\mu_{\cdot\wedge\tau_N}^{N+1})}_{\cd\wedge\tau_N}=(\sT(\bar\mu^{N+1})_{\cdot\wedge\tau_{N+1}})_{\cd\wedge\tau_N}=\sT_{N+1}(\bar\mu^{N+1})_{\cdot\wedge\tau_N}=\bar\mu_{\cdot\wedge\tau_N}^{N+1}.
$$
Thus $\bar\m^{N+1}_{\cdot\wedge\tau_N}$ is also a fixed point of $\sT_N$. By the uniqueness of the fixed point for $\sT_N$ we must have
$\bar\mu_{\cdot\wedge\tau_N}^{N+1}=\bar\mu^N$ on $[0,\t_N]$. That is, $\bar\m^{N+1}$ is an ``extension" of $\bar\m^N$.

The rest of the proof is now standard. We can ``patching" all the $\bar\m$'s together by defining a measure-valued process
$$\bar\mu_t:=\bar\mu_t^N, \ \ \ t\in[0,\tau_N], \ \ \  N\geq1.
$$
Then $\bar\m$ is well-defined on $[0,T]$ and one can easily check $\bar\m\in S_{\mathbb{F}^Y}^2(\sP_1)$. Furthermore, using
(\ref{15}) again we have, for any $N\geq1$,
\beaa
\sT(\bar\mu)_{\cdot\wedge\tau_N}=\sT(\bar\mu_{\cdot\wedge\tau_N})_{\cdot\wedge\tau_N}=\sT(\bar\mu^N)_{\cdot\wedge\tau_N}=\sT_N(\bar\mu^N)=\bar\mu^N=\bar\mu_{\cdot\wedge\tau_N}.
\eeaa
Thus, $\bar\m$ is a fixed point of $\sT$ on $[0,T]=\cup_{N=1}^\infty[0,\t_N]$.  Finally, note that if $\n$ is another fixed point of $\sT$, then by definition, for each $N\ge 1$, $\n_{\cd\wedge \t_N}$ must be
a fixed point of $\sT_N$. The uniqueness of the fixed point  then implies that $\bar\m_{\cd\wedge \t_N}=\n_{\cd\wedge \t_N}$,
which in turn  implies the uniqueness of the fixed point of $\sT$.
The proof is now complete.
\qed

Now let $\bar\m\in S_{\mathbb{F}^Y}^2(\sP_1)$ be the fixed point of $\sT$, and denote  $ (X,L):=(X^{\bar\mu},L^{\bar\mu})$. Recalling the
construction of $\sT(\bar\mu)$ we see that,  under $\hQ^0$, the couple of processes $(X, L)$ satisfies the following SDE:
\bea
\label{16}
\left\{\ba{lll}
dX_t=\sigma(t, X_{\cdot\wedge t}, Y_{\cdot\wedge t}, \bar\mu_{\cdot\wedge t} )dB_1^1,\qq & X_0=x_0,\\
dL_t=h(t, X_t, Y_{\cdot\wedge t})L_tdY_t, & L_0=1, \qq t\in[0,T].
\ea\right.
\eea
Furthermore, by construction (\ref{muth}) we see that $\bar\m_t=\hP\{X_t\in \cdot |\mathcal{F}_t^Y\}$, where $d\hP:=L_Td\hQ^0$. Now, denoting $\bar\m=\m^{X|Y}$,   we have the following theorem.
\begin{thm}
\label{exist}
Assume that the Assumptions \ref{assump1} and \ref{assump2} are in force. Then the SDE (\ref{SDE1}) possesses a weak solution.
\end{thm}

\noindent{\it Proof.} First note that, given the fixed point $\bar\m$ of the mapping $\sT$, and the corresponding solution $(X,L)$ to
SDE (\ref{16}) on the probability space $(\O, \cF, \hQ^0)$, if we define $d\hP:=L_Td\hQ^0$, and
$B_t^2:=Y_t-\int_0^th(s, X_s, Y_{\cdot\wedge s})ds$, $t\in[0,T]$, then the Girsanov theorem tells us that
the process $(B^1, B^2)$ is an $(\mathbb{F},\hP)$-Brownian motion.

Now, recalling that $\m^{X|Y}_t(\cd)=\bar\m_t(\cd)=\hP\{X_t\in \cdot |\mathcal{F}_t^Y\}$, $t\in[0,T]$,  we have, for $t\in[0,T]$,
\beaa
\left\{\ba{lll}
dX_t =\sigma(t, X_{\cdot\wedge t}, Y_{\cdot\wedge t}, \mu_{\cdot\wedge t}^{X|Y})dB_t^1,\qq & X_0=x_0;\\
dY_t =h(t, X_t, Y_{\cdot\wedge t})dt+dB_t^2, &Y_0=0.
\ea\right.
\eeaa
In other words, the six-tuple
$ (\O , \mathcal{F}, \mathbb{F}, \hP, (B^1, B^2), (X,Y))$ is a weak solution of SDE (\ref{SDE1}), proving the theorem.
\qed
%
%

\section{Uniqueness in law}
\setcounter{equation}{0}

In this section we shall address that last issue of the weak well-posedness of SDE (\ref{SDE1}). Namely, we shall prove that the weak solution of (\ref{SDE1}) is unique in law.
Our main idea extends the one in our previous work \cite{BLM} in a non-trivial way. That is, we note the fact that if $(X, Y, \hP)$ is a weak solution
to (\ref{SDE1}), and $\m^{X|Y}$ is the conditional law of $X$ given $\hF^Y$, under $\hP$, then
as we argued before we must have  $\m^{X|Y}\in S_{\mathbb{F}^{Y}}^2(\sP_1)\subset \hL^0_{\hF^Y}(\sC_T)$.
Therefore, there exists a progressively
measurable Borel functional $\Phi: \hC_T\rightarrow \sC_T(\sP_1)$, such that
\bea
\label{mui}
\m^{X |Y }_t=\Phi(Y )_t =\Phi(Y_{\cdot\wedge t} )_t ,\ \ t\in[0,T],\ \ \hP\mbox{-a.s.}
\eea
We shall use this function $\Phi$ as the bridge to connect two weak solutions, and then argue that they must be unique in law.
More precisely, we have the following theorem.
 \begin{thm} Assume the Assumptions \ref{assump1} and \ref{assump2}.
 Let
 $(\Omega^i, \mathcal{F}^i, \mathbb{F}^i, \hP^i, (B^{1,i},B^{2,i}), (X^i,Y^i)),\ i=1,2,$
 be two weak solutions of (\ref{SDE1}). Then, it holds that
 $$\hP^1 \circ (B^{1,1}, B^{2,1},X^1, Y^1)^{-1}=\hP^2 \circ (B^{1,2}, B^{2,2},X^2, Y^2)^{-1}.$$
 \end{thm}
\noindent
{\it{Proof.}}
Consider the following SDEs on $(\O^i, \hF^i,  \hP^i)$, $ i=1,2$, respectively:
\begin{equation}\label{17}
d\hat L_t^i=-\hat L_t^ih(t, X_t^i, Y_{\cdot\wedge t}^i)dB^{2, i}_t,\q \bar L_0^i=1, \q t\in[0,T].
\end{equation}
(Note the difference between this SDE and the one in (\ref{SDE2})!) Since $h$ is bounded, we know that $\hE^{\hP^i}[ \hat L_T^i ]=1$,  and
$d\hQ^i:=\hat L_T^i d\hP^i$ defines a probability measure such that $(B^{1,i},Y^i)$ is an $(\hF^i,\hQ^i)$-Brownian motion, $ i=1,2$.
Denote $L^i=[\hat L^i]^{-1}$.

Now let $\Phi^i: \hC_T\rightarrow \sC_T(\sP_1)$, $i=1,2$, be the progressively measurable Borel functionals, such that (\ref{mui})
holds for $ \m^{X^i|Y^i}$, $i=1,2$, respectively.
Then, the process $(X^1, L^1=[\hat L^1]^{-1})$ must satisfy the following SDE on $(\O^1, \cF^1, \hF^1, \hQ^1)$:
\begin{equation}
\label{18}
\left\{\ba{lll}
dX_t^1=\sigma(t, X_{\cdot\wedge t}^1, Y_{\cdot\wedge t}^1, \Phi_{\cdot\wedge t}^1(Y_{\cdot\wedge t}^1))dB_t^{1,1},\q & X_0^1=x_0,\\
dL_t^1=L_t^1h(t, X_t^1, Y_{\cdot\wedge t}^1)dY_t^1,  &L_0^1=1, \q  t\in[0,T].
\ea\right.
\end{equation}
Note that under $\hQ^1$, $(B^{1,1}, Y^1)$ is a Brownian motion, thus (\ref{18}) is just an SDE with random coefficients, and
under the Assumptions \ref{assump1} and \ref{assump2}, it has a pathwisely unique strong solution. That is, there exists a progressively measurable Borel functional $\psi: \hC^2_T\mapsto
\hC^2_T$, such that $(X^1,L^1) = \psi(B^{1,1},Y^1)$.

We now consider the following auxiliary SDE on the filtered space $(\O^2, \cF^2, \hF^2, \hQ^2)$:
\bea
\label{19}
\left\{\ba{lll}
d{\bar{X}}_t^2=\sigma(t, \bar{X}_{\cdot\wedge t}^2, Y_{\cdot\wedge t}^2, \Phi_{\cdot\wedge t}^1(Y_{\cdot\wedge t}^2))dB_t^{1,2},\qq & \bar{X}_0^2=x_0,\\
d{\bar{L}}_t^2=\bar{L}_t^2h(t,\bar{X}_t^2, Y_{\cdot\wedge t}^2)dY_t^2, & \bar{L}_0^2=1, \qq t\in[0,T].
\ea\right.
\eea
Note that SDE (\ref{19}) actually has the same coefficients as (\ref{18}), hence by pathwise uniqueness we deduce that
$(\bar{X}^2,\bar{L}^2)=\psi(B^{1,2},Y^2)$ as well.
But since $\hQ^1\circ (B^{1,1},Y^1)^{-1} = \hQ^2\circ (B^{1,2},Y^2)^{-1}$ is the Wiener measure on $(\O^0, \cF^0)=(\hC^2_T,\sB(\hC^2_T))$, we conclude that
\begin{equation}
\label{20}
\hQ^1\circ((B^{1,1},Y^1),(X^1,L^1))^{-1} = \hQ^2\circ((B^{1,2},Y^2),(\bar{X}^2,\bar{L}^2))^{-1}.
\end{equation}

Our next step is to use $(\bar X^2, \bar L^2)$ to build a bridge that links the laws of $(X^1, L^1)$ and $(X^2, L^2)$. To this end,
 let us now define a new probability $\bar{\hP}^2$ by $d\bar{\hP}^2 = \bar{L}^2_T d\hQ^2$, and consider the conditional
law $\bar\m^2=\bar\m^{\bar X^2|Y^2}$, under $\bar\hP^2$. That is, for  $ A\in\sB(\mathbb{R})$, it holds that
\bea
\label{bmu2}
\bar{\mu}^2_t(A)=\bar\m^{\bar X^2|Y^2}_t(A):=\bar{\hP}^2\big\{\bar{X}^2_t\in A\big|\mathcal{F}^{Y^2}_t\big\} = \frac{\hE^{\hQ^2}\big[\bar{L}^2_t\cdot\mathbf{1}_{\{\bar{X}^2_t\in A\}}\big|\mathcal{F}^{Y^2}_T\big]}{\hE^{\hQ^2}\big[\bar{L}^2_t\big|\mathcal{F}^{Y^2}_T\big]}.
\eea
We shall assume without loss of generality that $\bar\m^2$ is a regular conditional probability. As before, we can show that
  $\bar{\mu}^2\in S^2_{\mathbb{F}^Y}(\sP_1)$ (under $\hQ^2$). Furthermore, it holds that
\begin{equation}
\label{21}
\bar{\mu}^2_t(\cd) = \Phi^1_t(Y^2_{\cdot\wedge t})(\cd),\qquad t\in[0,T].
\end{equation}
Indeed, recall that $\Phi^1: \hC_T\rightarrow \sC_T(\sP_1)$ and observe that, for any bounded Borel functionals $\varphi:\hR\mapsto \hR$ and $ f:\hC_T\to \hR$, (\ref{20}) implies that
\bea
\label{bmu3}
&&\hE^{\bar{\hP}^2}\Big[\varphi\big(Y^2_{\cdot\wedge t}\big)\int_\mathbb{R} f(x) \Phi^1_t(Y^2_{\cdot\wedge t})(dx)\Big] = \hE^{\hQ^2}\Big[\bar{L}^2_t\varphi\big(Y^2_{\cdot\wedge t}\big)\int_\mathbb{R} f(x)  \Phi^1_t(Y^2_{\cdot\wedge t})(dx)\Big]
\nonumber\\
& =& \hE^{\hQ^1}\Big[L^1_t\varphi\big(Y^1_{\cdot\wedge t}\big)\int_\mathbb{R} f(x) \Phi^1_t(Y^1_{\cdot\wedge t})(dx)\Big].
\eea
Recalling that $\Phi^1_t(Y^1_{\cdot\wedge t})= \mu^{X^1|Y^1}_t(\cd)= \hP^1\big\{{X}^1_t\in \cd\,\big|\mathcal{F}^{Y^1}_t\big\}$,
$t\in[0,T]$, we have
$$ \int_\mathbb{R} f(x) \Phi^1_t(Y^1_{\cdot\wedge t})(dx)=\hE^{\hP^1}\big[f(X^1_t)\big|\mathcal{F}^{Y^1}_t\big].
$$
Thus (\ref{bmu3}) now reads
\bea
\label{EQ1}
&&\hE^{\bar{\hP}^2}\Big[\varphi\big(Y^2_{\cdot\wedge t}\big)\int_\mathbb{R} f(x) \Phi^1_t(Y^2_{\cdot\wedge t})(dx)\Big]
=\hE^{\hQ^1}\Big[L^1_t\varphi\big(Y^1_{\cdot\wedge t}\big)\hE^{\hP^1}\big[f(X^1_t)\big|\mathcal{F}^{Y^1}_t\big]\Big]\\
&=& \hE^{\hP^1}\big[\varphi\big(Y^1_{\cdot\wedge t}\big)\hE^{\hP^1}\big[f(X^1_t)\big|\mathcal{F}^{Y^1}_t\big]\big]
= \hE^{\hP^1}\big[\varphi\big(Y^1_{\cdot\wedge t}\big)f(X^1_t)\big]
= \hE^{\hQ^1}\Big[L^1_t\varphi\big(Y^1_{\cdot\wedge t}\big)f(X^1_t)\Big].\nonumber
\eea
On the other hand, (\ref{EQ1}), together with (\ref{20}), also shows that
\beaa
\label{EQ2}
&&\hE^{\bar{\hP}^2}\Big[\varphi\big(Y^2_{\cdot\wedge t}\big)\int_\mathbb{R} f(x)  \Phi^1_t(Y^2_{\cdot\wedge t})(dx)\Big]
= \hE^{\hQ^1}\big[L^1_t\varphi\big(Y^1_{\cdot\wedge t}\big)f(X^1_t)\big]\\
&=& \hE^{\hQ^2}\big[\bar{L}^2_t\varphi\big(Y^2_{\cdot\wedge t}\big)f(\bar{X}^2_t)\big]= \hE^{\bar{\hP}^2}\big[\varphi\big(Y^2_{\cdot\wedge t}\big)f(\bar{X}^2_t)\big]\\
&=& \hE^{\bar{\hP}^2}\Big[\varphi\big(Y^2_{\cdot\wedge t}\big)\hE^{\bar{\hP}^2}\big[f(\bar{X}^2_t)\big|\mathcal{F}^{Y^2}_t\big]\Big]
= \hE^{\bar{\hP}^2}\Big[\varphi\big(Y^2_{\cdot\wedge t}\big)\int_\mathbb{R} f(x) \bar{\mu}^2_t(dx)\Big]. \nonumber
\eeaa
Since both $\f$ and $f$ are arbitrary, we have proved the claim (\ref{21}).
We note that using (\ref{21}) SDE (\ref{19})  can be rewritten as
\bea
\label{22}
\left\{ \ba{lll}
         d\bar{X}^2_t = \sigma(t,\bar{X}^2_{\cdot\wedge t},Y^2_{\cdot\wedge t},\bar{\mu}^2_{\cdot\wedge t})dB^{1,2}_t,\qq & \bar{X}^2_0 = x_0,  \\
         d\bar{L}^2_t = \bar{L}^2_t h(t,\bar{X}^2_t,Y^2_{\cdot\wedge t})dY^2_t, &\bar{L}^2_0 = 1, \qquad t\in[0,T],
 \ea\right.
\eea
with $\bar{\mu}^2_{\cdot\wedge t}(\cd) = \bar{\hP}^2\big\{\bar{X}^2_t\in \cdot\big|\mathcal{F}^{Y^2}_t\big\},\ t\in[0,T]$, which
satisfies (\ref{bmu2}).

Our final observation is that, by the construction of the solution mapping $\sT$ (\ref{muth}) and the definition of $\m^2$, we see that
both $\m^2$ and $\bar\m^2$ are in
$S^2_{\hF^{Y^2}}(\sP_1)$ under the probability $\hQ^2$,  and they  satisfy
$$\sT(\bar{\mu}^2)=\bar{\mu}^2,\qq \sT({\mu}^2)={\mu}^2.$$
Namely, both $\m^2$ and $\bar \m^2$ are the fixed points of the solution mapping $\sT$. Thus, the uniqueness of the fixed point
implies that $\bar{\mu}^2=\mu^2$.

Finally, recall that the process $(X^2, L^2=[\hat L^2]^{-1})$ satisfies the SDE
\begin{equation}
\label{23}
\left\{ \ba{lll}
         d{X}^2_t = \sigma(t,{X}^2_{\cdot\wedge t},Y^2_{\cdot\wedge t},{\mu}^2_{\cdot\wedge t})dB^{1,2}_t, \qq &{X}^2_0 = x_0,  \\
         d{L}^2_t = {L}^2_t h(t,{X}^2_t,Y^2_{\cdot\wedge t})dY^2_t, & {L}^2_0 = 1, \qquad t\in[0,T],
         \ea \right.
\end{equation}
where  ${\mu}^2_{\cdot\wedge t} = {\hP}^2\big\{{X}^2_t\in\cdot\big|\mathcal{F}^{Y^2}_t\big\}$, $t\in[0,T]$,
$d{\hP}^2 = {L}^2_T d\hQ^2$.
Consequently, both SDEs (\ref{22}) and (\ref{23}) are defined on $(\O^2, \cF^2, \hQ^2)$, have the same coefficients (given $\m^2$),
and are driven by the same $(\hF^2, \hQ^2)$-Brownian motion $(B^{1,2},Y^2)$. Thus the pathwise uniqueness of SDE (given $\m^2$) leads to that
$(\bar{X}^2,\bar{L}^2)\equiv (X^2,L^2)$, $\hQ^2$-a.s.

Consequently, we now have $d\bar\hP^2 = \bar{L}^2_T d\hQ^2 = {L}^2_T d\hQ^2 = d{\hP}^2$.
Combining this with (\ref{20}) we get
\begin{equation}\label{24}
\hQ^1\circ(B^{1,1},Y^1,X^1,L^1)^{-1}
= \hQ^2\circ(B^{1,2},Y^2,\bar{X}^2,\bar{L}^2)^{-1}
= \hQ^2\circ(B^{1,2},Y^2,X^2,L^2)^{-1}.
\end{equation}
Since  $\displaystyle{B^{2,1}_t = Y^1_t - \int^t_0 h(s,X^1_s,Y^1_{\cdot\wedge s})ds}$,\ $\displaystyle{B^{2,2}_t = Y^2_t - \int^t_0 h(s,X^2_s,Y^2_{\cdot\wedge s})ds},\ t\in[0,T] $, we obtain from (\ref{24}) that
\beaa
\hP^1\circ(B^{1,1},B^{2,1},X^1,Y^1)^{-1}
= \hP^2\circ(B^{1,2},B^{2,2},X^2,Y^2)^{-1}.
\eeaa
The proof is now complete.
\qed


\begin{thebibliography}{1}
\bibitem{BLM}
Buckdahn, R., Li, J.,  and Ma, J.,  {\it A Mean-field Stochastic Control Problem with   Partial  Observations},
Ann. Appl. Probab., {\bf 27} (2017), no. 5, 3201-3245.


\bibitem{Billing}
Billingsley, P., {\sl Convergence of Probability Measures}, John Wiley \& Sons, (2013).

\bibitem{BCEH}
Briand, P.,  Cardaliaguet, P., \'{E}ric Chaudru de Raynal, P., and Hu, Y., {\it Forward and Backward Stochastic
Differential Equations with normal constraint in law}, Stoch. Proc. Appl., {\bf 130} (2020), no. 12, 7021-7097.


\bibitem{CDL}
Carmona, R., Delarue, F., and Lacker, D., {\it Mean-field games with common noise,} Ann. Probab., {\bf 44} (2016), no. 6, 3740-3803.

\bibitem{CG}
Coghi, M. and Gess, B., {\it Stochastic nonlinear Fokker-Planck equations},
Nonlinear Analysis, {\bf 187} (2019), no.  259-278.

\bibitem{Ed}
Edwards, D.~A., {\it On the Kantorovich-Rubinstein Theorem}, Expositiones Mathematicae, {\bf 29} (2011), 387--398.

\bibitem{EK}
Ethier, S., and Kurtz,  T.~G., {\sl Markov Processes, Characterization and Convergence}, John Williams \& Sons Inc.,
(1986).

\bibitem{Flandoli}
Flandoli, F., {\sl Compact Sets in the Space of Measure-valued Functions}, Lecture notes??

\bibitem{HSS}
Hammersley, W.~R.,  Siska, D., and Szpruch, L.,
{\it Weak existence and uniqueness for Mckean-Vlasov SDEs with common
  noise}, Ann. Probab., {\bf49} (2021), no. 2, 527-555.

\bibitem{KA}
Kantorovich, L.~V. and Akilov, G.~P., {\sl Functional Analysis}, 2nd ed., Pergamon Press, Oxford, 1982.

\bibitem{KR}
Kantorovich, L.~V. and Rubinstein, G., {\it On a space of completely additive functions}, Vestnik Leningradskogo Universiteta, {\bf 13} (1958), no. 7,
52-59.

\bibitem{Krylov}
Krylov, N.~V., {\it An analytic approach to SPDEs}. {\sl Stochastic partial differential equations: six perspectives}, 185-242, Math. Surveys Monogr., 64, Amer. Math. Soc., Providence, RI, 1999.

\bibitem{LS}
Ledger, S. J. and S{o\neg\neg\neg/}jmark, A., {\it At the Mercy of the Common Noise: Blow-ups in a Conditional McKean-Vlasov Problem}, Electron. J. Probab., (2021), no. 26, 1-39.

\bibitem{MSZ}
Ma, J.,  Sun, R., and  Zhou, Y., {\it Kyle-Back Equilibrium Models and Linear Conditional Mean-field SDEs},
SIAM J. Control Optim., {\bf 56} (2018), no. 2, 1154 -1180.\end{thebibliography}
\end{document}